\documentclass{amsart}
\usepackage{graphics,amssymb,amscd,amsthm,verbatim}
\input epsf                 

\newtheorem{proposition}{Proposition}[section]
\newtheorem{theorem}[proposition]{Theorem}

\newtheorem{lemma}[proposition]{Lemma}

\newtheorem{prop}[proposition]{Proposition}

\newtheorem{conj}[proposition]{Conjecture}

\newcommand{\reals}{\mathbb R}

\newcommand{\A}{{\mathcal A}}

\newcommand{\C}{{\mathcal C}}

\newcommand{\F}{{\mathcal F}}
\newcommand{\Po}{{\mathcal P}}

\newcommand{\covers}{{\,\,\,\cdot\!\!\!\! >\,\,}}
\newcommand{\covered}{{\,\,<\!\!\!\!\cdot\,\,\,}}

\newcommand{\set}[1]{{\left\lbrace #1 \right\rbrace}}

\newcommand{\join}{\vee}
\newcommand{\meet}{\wedge}
\newcommand{\Cg}{\mathrm{Cg}}
\newcommand{\tw}{\mathrm{tw}}
\newcommand{\Con}{\mathrm{Con}}

\newcommand{\Irr}{\mathrm{Irr}}

\newcommand{\st}{\mathrm{st}}
\newcommand{\slo}{\mathrm{sl}}

\newcommand{\pidown}{\pi_\downarrow}
\newcommand{\piup}{\pi^\uparrow}

\newcommand{\case}[2]{\item[Case #1:] #2 \\ }

\newcommand{\mapname}[1]{\stackrel{#1}{\longrightarrow}}

\newlength{\lsmash} % this sets how much the 2nd overline/underline in \up/\down is mashed down onto the first
\newlength{\myheight}
\newlength{\mydepth}

\newcommand{\up}[1]{
\settoheight{\myheight}{\ensuremath{\overline{#1}}}
\addtolength{\myheight}{-\lsmash}
\overline{
\protect \raisebox{0 pt}[\myheight][0 pt]{\ensuremath{\overline{#1}}}
}
}

\newcommand{\down}[1]{
\settodepth{\mydepth}{\ensuremath{\underline{#1}}}
\addtolength{\mydepth}{-\lsmash}
\underline{
\protect \raisebox{0 pt}[0 pt][\mydepth]{\ensuremath{\underline{#1}}}
}
}

\newcommand{\upwide}[1]{
\settoheight{\myheight}{\ensuremath{\overline{\,#1\,}}}
\addtolength{\myheight}{-\lsmash}
\overline{
\protect \raisebox{0 pt}[\myheight][0 pt]{\ensuremath{\overline{\,#1\,}}}
}
}

\newcommand{\downwide}[1]{
\settodepth{\mydepth}{\ensuremath{\underline{\,#1\,}}}
\addtolength{\mydepth}{-\lsmash}
\underline{
\protect \raisebox{0 pt}[0 pt][\mydepth]{\ensuremath{\underline{\,#1\,}}}
}
}

\newcommand{\Pge}{{\Phi_{\ge -1}}}
\newcommand{\ck}{^\vee}
\newcommand{\G}{{\vec{G}}}
\newcommand{\ep}{\varepsilon}

\newcommand{\br}[1]{\langle #1 \rangle}

\begin{document}
\title{Cambrian Lattices}

\author{Nathan Reading}
\address{
Mathematics Department\\
       University of Michigan\\
       Ann Arbor, MI 48109-1109\\
USA}
\thanks{The author was partially supported by NSF grant DMS-0202430.}
\email{nreading@umich.edu}
\urladdr{http://www.math.lsa.umich.edu/$\sim$nreading/}
%\keywords{NEED THESE}

\begin{abstract}
For an arbitrary finite Coxeter group $W,$\ we define the family of Cambrian lattices for $W$ as quotients of the weak order on $W$ 
with respect to certain lattice congruences.
We associate to each Cambrian lattice a complete fan, which we conjecture is the normal fan of a polytope combinatorially isomorphic to the
generalized associahedron for $W.$\
In types~A and~B we obtain, by means of a fiber-polytope construction, combinatorial realizations of the Cambrian lattices in terms of 
triangulations and in terms of permutations.
Using this combinatorial information, we prove in types~A and~B that the Cambrian fans are combinatorially isomorphic to the normal fans of 
the generalized associahedra and that one of the Cambrian fans is linearly isomorphic to Fomin and Zelevinsky's construction of the normal
fan as a ``cluster fan.'' 
Our construction does not require a crystallographic Coxeter group and therefore suggests a definition, at least on the level of cellular spheres, of a 
generalized associahedron for any finite Coxeter group.
The Tamari lattice is one of the Cambrian lattices of type~A, and two ``Tamari'' lattices in type~B are identified and characterized in terms of 
signed pattern avoidance.
We also show that open intervals in Cambrian lattices are either contractible or homotopy equivalent to spheres.
\end{abstract}

\maketitle

\section{Introduction}
\label{intro}

The associahedron or Stasheff polytope is a polytope whose vertices are counted by the Catalan numbers.
The Tamari lattice is a lattice whose Hasse diagram is the 1-skeleton of the associahedron.
There is a $W$-Catalan number~\cite[Remark 2]{NonCrossing} associated to any Coxeter group $W$ and the ordinary Catalan numbers 
are associated to irreducible Coxeter groups of type~A (the symmetric groups).
Chapoton, Fomin and Zelevinsky~\cite{gaPoly,ga} have recently generalized the associahedron to all finite crystallographic Coxeter groups.
Earlier, the \mbox{type-B} associahedron, or cyclohedron, was defined by Bott and Taubes~\cite{B-T} and Simion~\cite{B assoc}.

Simion~\cite[\S 5]{B assoc} asked if the vertices of the \mbox{type-B} associahedron could be partially ordered so that the Hasse diagram of the 
partial order is the 1-skeleton of the \mbox{type-B} associahedron.
Reiner~\cite[Remark 5.4]{Equivariant} used an equivariant fiber polytope construction to identify a family of maps from the \mbox{type-B} 
permutohedron to the \mbox{type-B} associahedron, in analogy to well-known maps in type~A.
He also asked whether one of these maps can be used to define a partial order on the vertices of the \mbox{type-B} associahedron with similar properties
to the Tamari lattice, such that the map from the \mbox{type-B} permutohedron to the \mbox{type-B} associahedron shared the pleasant properties of the 
corresponding map in type~A (see~\cite{Nonpure II}).

The starting point of the present research is the observation that the Tamari lattice is a lattice-homomorphic image of the right weak order on 
the symmetric group.
This fact has, to our knowledge, never appeared in the literature, although essentially all the ingredients of a proof were assembled by Bj\"{o}rner
and Wachs in~\cite{Nonpure II}.
This lattice-theoretic point of view suggests a search among Reiner's maps to determine which induces a lattice homomorphism on the weak 
order.
Surprisingly, for each of these maps, one can choose a vertex to label as the identity element such that the map induces a lattice homomorphism, 
and thus each defines a lattice structure on the \mbox{type-B} associahedron.
Furthermore, the analogous family of maps in type~A yields a family of lattices on vertices of the \mbox{type-A} associahedron.
A close look at the lattice homomorphisms in types~A and~B leads to a type-free generalization of these families of lattices which we call 
{\em Cambrian lattices}.\footnote{
The Cambrian layer of rocks marks a dramatic increase in the diversity in the fossil record and thus the sudden profusion of Catalan-related 
lattices arising from the single (Pre-Cambrian) example of the Tamari lattice might fittingly be called ``Cambrian.''}

Let $G$ be the Coxeter diagram of a finite Coxeter group $W.$\
An {\em orientation} $\G$ of $G$ is a directed graph with the same vertex set as $G$, with one directed edge for each edge of $G$.
For each orientation $\G$ of $G$, there is a Cambrian lattice which we now proceed to define.
For a directed edge $s\to t$ in $\G$ with label $m(s,t)$, require that $t$ be equivalent to the element of $W$ represented by the word 
$tsts\cdots$ with $m(s,t)-1$ letters.
The {\em Cambrian congruence} $\Theta(\G)$ associated to $\G$ is the smallest lattice congruence of the right weak order on $W$ satisfying this 
requirement for each directed edge in $\G$.

\begin{figure}[ht]
\centerline{\scalebox{.8}{\epsfbox{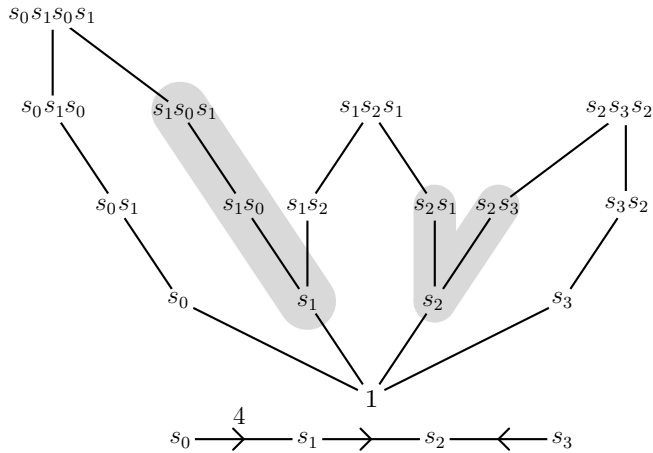}}}   
\caption{How an orientation specifies equivalences.}
\label{B4cong}
\end{figure}

\begin{figure}[ht]
\centerline{\scalebox{1}{\epsfbox{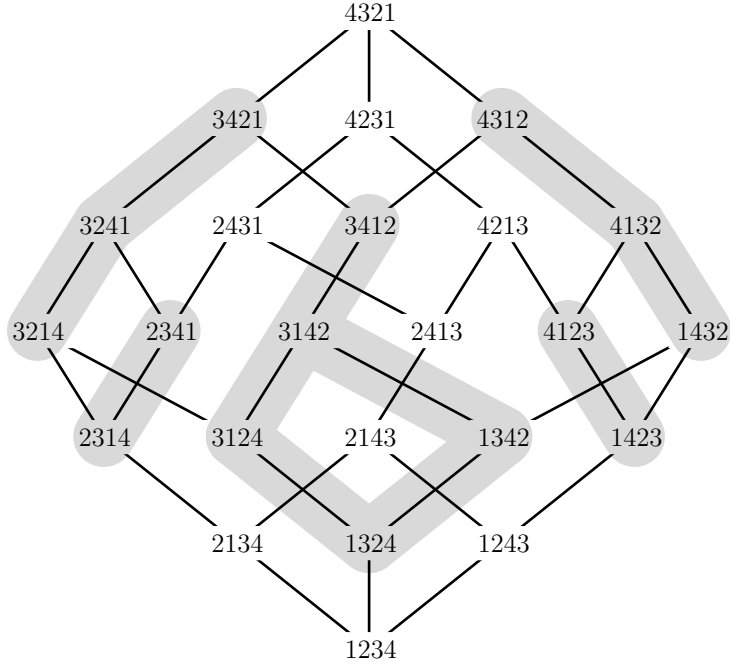}}}   
\caption{A Cambrian congruence on $S_4$.}
\label{A3camb_cong}
\end{figure}

Figure~\ref{B4cong} shows part of the weak order on $B_4$, with shaded edges indicating the required equivalences for the orientation shown.
Figure~\ref{A3camb_cong} shows the Cambrian congruence of the weak order on $S_4$ arising from the orientation with arrows $2134\to1324$ and $1243\to1324$ (cf. Figure~\ref{A3tri}).
This is the smallest lattice congruence with $1324\equiv3124$ and $1324\equiv1342$.
Each shaded edge in the Hasse diagram indicates a covering pair $x\covered y$ such that $x\equiv y$ and the congruence in question is the transitive closure of these equivalences of covering pairs.

The Cambrian lattice $\C(\G)$ is defined to be the right weak order on $W$ modulo the congruence $\Theta(\G)$.
Figure~\ref{A3camb_tri} in Section~\ref{Cam A} shows the Cambrian lattice arising from the Cambrian congruence of Figure~\ref{A3camb_cong}.
Two Cambrian lattices are isomorphic (respectively anti-isomorphic) exactly when the associated directed diagrams are isomorphic (respectively 
anti-isomorphic), taking edge labelings into account (Theorem~\ref{iso}).

For any finite Coxeter group $W,$\ let $\F$ be the fan defined by the reflecting hyperplanes of $W,$\ and identify the maximal cones with the 
elements of $W.$\
In~\cite{con_app} a fan $\F_\Theta$ is defined for any lattice congruence $\Theta$ of the weak order on $W,$\ such that the maximal cones 
of $\F_\Theta$ are the unions over congruence classes of the maximal cones of $\F$.
Thus the maximal cones of $\F_\Theta$ correspond to elements of the lattice quotient $W/\Theta$.
Let $\F(\G)$ denote the {\em Cambrian fan} $\F_{\Theta(\G)}$.
\begin{conj}
\label{fan conj}
For any finite Coxeter group $W$ and any orientation $\G$ of the associated Coxeter diagram, the fan $\F(\G)$ is the normal fan of a 
simple convex polytope.
If $W$ is crystallographic then this polytope is combinatorially isomorphic to the generalized associahedron for $W.$\
\end{conj}
Each statement in the following conjecture is weaker than Conjecture~\ref{fan conj}, and proofs of any of these  weaker conjectures 
would be interesting.
\begin{conj}
\label{weaker fan conjs}
Let $W$ be a finite Coxeter group with digram $G$:
\begin{enumerate}
\item[a. ]If $W$ is crystallographic then given any orientation $\G$ of $G$, the fan $\F(\G)$ is combinatorially isomorphic to the normal fan of the generalized associahedron 
for $W.$\
\item[b. ]Given any orientation $\G$ of $G$, the fan $\F(\G)$ is combinatorially isomorphic to the normal fan of some polytope.
\item[c. ]Cambrian fans arising from different orientations of $G$ are combinatorially isomorphic.
\item[d. ]Cambrian fans arising from different orientations of $G$ have the same number of maximal cones.
\item[e. ] Given any orientation $\G$ of $G$, the fan $\F(\G)$ is simplicial.
\end{enumerate}
\end{conj}

\begin{figure}[ht]
\centerline{\epsfbox{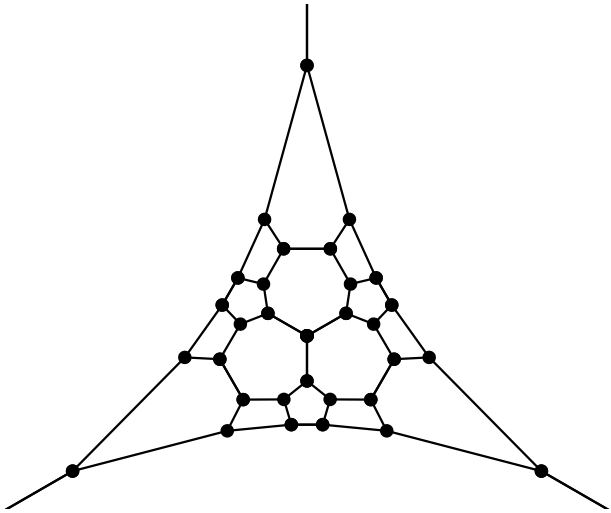}}   
\caption{The $1$-skeleton of the $H_3$-associahedron.  
The vertex at infinity completes the three unbounded regions to heptagons.}
\label{H3}
\end{figure}

The fan $\F(\G)$ has an associated regular cellular sphere and a dual regular cellular sphere $\Gamma(\G)$.
If Conjecture~\ref{fan conj} holds, then in particular, $\Gamma(\G)$ is a polytope, and the Cambrian fans offer an alternate definition
of the generalized associahedra.
This construction would lift the restriction to crystallographic Coxeter groups imposed by the definition in~\cite{ga}, giving the first definition\footnote{Since the initial version of this paper, it has become apparent that the definition from~\cite{ga} can be extended to the noncrystallographic case.  See~\cite[Section 5.3]{rsga}.}
of associahedra of types H and~I.
One can verify that Conjectures~\ref{weaker fan conjs}.c, \ref{weaker fan conjs}.d and~\ref{weaker fan conjs}.e hold for $H_3$, and they hold 
trivially for $I_2(m)$.
The associahedra for $H_3$ and $I_2(m)$ constructed in this way have the numbers of faces of each dimension one would expect from
generalized associahedra, and the facets of the $H_3$-associahedron are the correct generalized associahedra of lower dimension.
The $H_3$-associahedron has a dihedral combinatorial symmetry group with 12 elements, as one would expect based on the crystallographic 
cases (see~\cite[Theorem 1.2]{ga}).
The $I_2(m)$-associahedron is an $(m+2)$-gon and the 1-skeleton of the $H_3$-associahedron is pictured in Figure~\ref{H3}.

The Cambrian lattices and fans have the following properties which follow from the results of~\cite{con_app}.
First, $\C(\G)$ is a partial order induced on the maximal cones of $\F(\G)$ by a linear functional, and the Hasse diagram of 
$\C(\G)$ is isomorphic to the 1-skeleton of the dual sphere $\Gamma(\G)$.
Any linear extension of $\C(\G)$ is a shelling of the sphere associated to $\F(G)$.
The set of cones containing a given face $F$ of $\F(\G)$ is a closed interval in $\C(\G)$ called a {\em facial interval}.
Open intervals associated to non-facial closed intervals in $\C(\G)$ are contractible and open intervals associated to facial intervals are 
homotopy-equivalent to spheres with the dimension of the sphere 
depending on the dimension of the corresponding face of $\F(\G)$.
This in particular determined the M\"{o}bius function of all intervals in $\C(\G)$.
Furthermore, intervals are facial if and only if they are atomic (meaning that the join of the atoms of the interval is the top of the interval).
The face lattice of $\Gamma(\G)$ is isomorphic to the set of atomic intervals, partially ordered by inclusion.
All combinatorial information about $\Gamma(\G)$ is thus encoded in the Cambrian lattice $\F(\G)$.
If $\F(\G)$ is indeed simplicial then the corresponding simplicial sphere is flag.
If $\F(\G)$ is indeed polytopal, then since it refines the normal fan of $\F$, the associated polytope is a Minkowski summand of the 
$W$-permutohedron.

Because disconnected Coxeter diagrams lead to direct product decompositions of all of the relevant objects, it is enough to prove 
Conjecture~\ref{fan conj} in the case of connected Coxeter diagrams, or equivalently irreducible Coxeter groups.
We use (equivariant) fiber polytope constructions to work out the combinatorics of Cambrian lattices of types~A and~B in detail, both in terms of
triangulations and in terms of permutations (Theorems~\ref{tri A}, \ref{per A}, \ref{tri B} and~\ref{per B}).
In particular:
\begin{theorem}
\label{ab weak fan}
Conjecture~\ref{weaker fan conjs}.a holds when $W$ is of type~A or~B\@.
\end{theorem}

Fomin's and Zelevinsky's definition~\cite{ga} of generalized associahedra is in terms of clusters of roots in the root system
associated to a crystallographic Coxeter group $W.$\
The clusters arise, via certain piecewise linear maps $\tau_+$ and $\tau_-$, from a bipartition $G=I_+\cup I_-$ of the Coxeter diagram.
The clusters are the vertices of the $W$-associahedron, and the edges are pairs of clusters which differ by exchanging one root.
The cluster fan is the fan whose maximal cones are the cones generated by the clusters.
It is also the normal fan of the $W$-associahedron.
Using $\tau_+$ and $\tau_-$ to compare two clusters that form an edge, in Section~\ref{cluster} we define a partial order on the clusters 
called the {\em cluster poset}.
Naturally associated to the bipartition $G=I_+\cup I_-$ is an orientation of $G$ which we denote $I_+\longrightarrow I_-$, and call a 
{\em bipartite orientation}.
Specifically, any edge in $G$ connects an element $s$ of $I_+$ to an element $t$ of $I_-$, and we direct the edge $s\to t$.
Let $\alpha_i$ denote a simple root in the associated root system $\Phi$ and let $\omega_i$ be the corresponding fundamental weight.
For any $i\in I$, set $\ep(i)=+1$ if $i\in I_+$ and $\ep(i)=-1$ if $i\in I_-$.

\begin{conj}
\label{cluster conj}
Let $W$ be a finite crystallographic Coxeter group.
Then the linear map $\alpha_i\mapsto\ep(i)\omega_i$ maps the cluster fan for~$W$ isomorphically to the Cambrian fan for the orientation 
$I_+\longrightarrow I_-$ of the diagram of $W.$\
The Cambrian lattice for the same orientation is the cluster poset.
\end{conj}
This conjecture would in particular imply that the cluster poset is a lattice, that it is induced on the vertices of the generalized associahedron by
a linear functional, that its Hasse diagram is isomorphic to the 1-skeleton of the generalized associahedron, and that it has the pleasant 
topological and combinatorial properties of the Cambrian lattices, as described above.
General proofs of any of these weaker statements would also be interesting.

Conjecture~\ref{cluster conj} can be proven in types~A and~B (Section~\ref{fans}).
This provides a proof of Conjecture~\ref{fan conj} in the special case where $\G$ is a bipartite orientation of the diagram of a Coxeter group of 
type~A or~B\@.
As further support for Conjecture~\ref{cluster conj}, we prove the following fact which would be a consequence of Conjecture~\ref{cluster conj}.
\begin{theorem}
\label{cluster refine}
The linear map $\alpha_i\mapsto\ep(i)\omega_i$ maps the cluster fan of a crystallographic Coxeter group $W$ to a fan that is refined by the normal fan of the $W$-permutohedron.
\end{theorem}
Marsh, Reineke and Zelevinsky~\cite{quiver} associate a fan to any Dynkin quiver such that the alternating quiver
is associated to the cluster fan.
Since their construction in~\cite{quiver} starts with the same combinatorial data as ours, it would be interesting to know if their fans coincide
with the Cambrian fans.

The combinatorial description of Cambrian lattices of type~A given in Section~\ref{Cam A} shows that the Tamari lattice is the 
\mbox{type-A} Cambrian lattice associated to a path directed linearly, that is, with the arrows all pointing the same direction.
Call this the {\em Tamari orientation} of the diagram.
By Theorem~\ref{iso} and the symmetry of the diagram, we recover the fact that the Tamari lattice is self-dual.
The Coxeter diagram for type~B is a path as well, but has an asymmetric edge-labeling.
There are two linear orientations of the \mbox{type-B} diagram, yielding two ``Type-B Tamari lattices,'' which are not isomorphic but dual to each other.
In support of their claim to the title of ``Tamari'' is the fact that each can be constructed as the sublattice of the weak order consisting of 
 signed permutations avoiding certain signed patterns (Theorem~\ref{B Tamari}).
This is analogous to the well-known realization of the \mbox{type-A} Tamari lattice in terms of pattern-avoidance~\cite{Nonpure II}.
Because the \mbox{type-B} Tamari elements are counted by the \mbox{type-B} Catalan numbers, this result has some bearing on a question posed by Simion in 
the introduction to~\cite{signed}, which asked for signed permutation analogues of counting formulas for restricted permutations.
Thomas~\cite{thomas}, working independently and roughly simultaneously, constructed the \mbox{type-B} Tamari lattice and proposed 
a type-D Tamari lattice.

Stasheff and Schnider~\cite{Sta95} gave a realization of the \mbox{type-A} associahedron by specifying facet hyperplanes, and Loday~\cite{Lo} 
determined the vertices of this realization.
The Cambrian fan for the \mbox{type-A} Tamari orientation is the normal fan of this realization of the associahedron (Section~\ref{fans}), 
and this observation proves Conjecture~\ref{fan conj} for the Tamari orientation in type~A\@.
Thus in type~A the Cambrian fans interpolate between the cluster fan and the normal fan of Stasheff's and Shnider's realization of the 
associahedron.

In general, the quotient of a lattice~$L$ with respect to some congruence is isomorphic to an induced subposet of~$L$, but need not be a sublattice.
However, in types~A and~B, the Cambrian lattices are sublattices of the weak order (Theorems~\ref{sub A} and~\ref{sub B}).
This fact was proven for the Tamari lattice in~\cite{Nonpure II}.
\begin{conj}
\label{sublattice conj}
For any finite Coxeter group $W$ and any orientation $\G$ of the associated Coxeter diagram, the Cambrian lattice $\C(\G)$ is a sublattice of 
the weak order on $W.$\
\end{conj}

The Cambrian lattices inherit any lattice property from the weak order which is preserved by homomorphisms.
Notably, the Cambrian lattices are congruence uniform, generalizing a theorem of Geyer~\cite{Geyer} on the Tamari lattice.

The left descent map $\delta:W\to 2^S$ is a lattice homomorphism from the right weak order on $W$ onto a boolean 
algebra~\cite{Barbut} (see also~\cite{congruence}).
Let $\eta$ be the canonical homomorphism from the weak order on $W$ to a Cambrian lattice.
Then $\delta$ factors through $\eta$ in the sense that there is a lattice homomorphism $\delta'$ from the Cambrian lattice to $2^S$ such that 
$\delta'\circ\eta=\delta:W\to 2^S$.
In types~A and~B we identify this map on triangulations (Theorems~\ref{A descents} and~\ref{B descents}).

This paper is the third in a series of papers beginning with~\cite{congruence} and~\cite{con_app}.
Each paper relies on the results of the preceding papers and cites later papers only for motivation or in the context of examples.
In a subsequent paper~\cite{sortable}, motivated by the study of Cambrian lattices, we define sortable elements of a Coxeter group~$W$ and use them to make a bijective connection between clusters and noncrossing partitions.
The connection between sortable elements and Cambrian lattices is not explicitly made in~\cite{sortable} but is the subject of~\cite{sort_camb}.

The remainder of the paper is organized as follows:
In Sections~\ref{weak} and~\ref{cong} we give background information on Coxeter groups, the weak order and congruences on the weak order.
In Sections~\ref{A assoc} and~\ref{eta A} we review the iterated fiber polytope construction mentioned above, and carefully analyze the 
resulting map from the permutohedron to the associahedron.
This analysis leads to the proof, in Section~\ref{Cam A}, of the \mbox{type-A} results described above.
In Section~\ref{Cam B} we review the equivariant fiber polytope construction and prove \mbox{type-B} versions of the results of Section~\ref{Cam A}.
Section~\ref{cluster} is a short description of clusters and a proof of Theorem~\ref{cluster refine}.
The proof of Conjecture~\ref{cluster conj} for types~A and~B is in Section~\ref{fans}, accompanied by an exact description of the Cambrian fan 
associated to the \mbox{type-A} Tamari lattice.

\section{The weak order on a Coxeter group}
\label{weak}
In this section we provide background information on Coxeter groups and the weak order, including a repetition in greater depth of a few 
concepts mentioned in the introduction.
For more information on Coxeter groups and weak order, see~\cite{orderings,Bourbaki,Humphreys}.

A {\em Coxeter system} $(W,S)$ is a group $W$ which can be presented as the group generated by the elements of the set $S$, subject to the 
relations $s^2=1$ for all $s\in S$ as well as the braid relations: for each $s,t\in S$ with $s\neq t$ there is a {\em pairwise order} $m(s,t)$ with 
$2\le m(s,t)\le\infty$ which is the smallest integer such that $(st)^{m(s,t)}=1$.
Typically, $W$ is called a {\em Coxeter group} with the understanding that there is some~$S$ for which $(W,S)$ is a Coxeter system.
A Coxeter system is encoded by a {\em Coxeter diagram} $G$, a graph whose vertex set is $S$, with an edge labeled $m(s,t)$ connecting~$s$ 
to~$t$ whenever $m(s,t)>2$.
Since $3$ is a common edge-label in Coxeter diagrams, by convention unlabeled edges are interpreted to mean $m(s,t)=3$.
A Coxeter group $W$ is {\em irreducible} if $G$ is connected.  
If $W$ is reducible, then it is isomorphic to the direct product of the irreducible Coxeter groups associated to the connected components of its 
diagram.
The elements of $S$ are called {\em simple reflections} and conjugates of simple reflections are called {\em reflections}.
Any element $w$ of $W$ can be written, non-uniquely, as a word in the alphabet $S$.
A minimal-length word for $w$ is called a {\em reduced word} for $w$, and the length of a reduced word is called the {\em length} of $w$ and
denoted $l(w)$.
The {\em left inversion set} of $w$, denoted $I(w)$ is the set of all (not necessarily simple) reflections $t$ such that $l(tw)<l(w)$.
The length $l(w)$ is equal to $|I(w)|$.
The {\em left descent set} of $w$ is $I(w)\cap S$.
There is an analogous {\em right inversion set} and {\em right descent set}.
The right descents of $w$ are the elements $s\in S$ such that $l(ws)<l(w)$.
If $W$ is a finite Coxeter group, there is a unique element $w_0$ of maximal length.

The {\em right weak order} on $W$ is the partial order which sets $v\le w$ if and only if $I(v)\subseteq I(w)$.
Alternately, the right weak order is the transitive closure of the cover relations $w\covered ws$ for every right descent $s$ of $w$.
There is a left weak order, isomorphic to the right weak order, but throughout this paper the term ``weak order'' always refers to the right weak 
order. 
For any $W,$\ the weak order is a meet-semilattice and if $W$ is finite then the weak order is a lattice, in which case $w_0$ is the unique
maximal element.

An element $\gamma$ of a finite lattice $L$ is join-irreducible if for any $X\subseteq L$ with $\gamma=\join X$, we have $\gamma\in X$.
Equivalently, $\gamma$ is join-irreducible if it covers exactly one element $\gamma_*$ of~$L$.
Thus an element $\gamma$ of a finite Coxeter group $W$ is join-irreducible in the weak order if and only if it has a unique right descent~$s$.
Given a join-irreducible $\gamma$ with right descent $s$, the {\em associated left reflection} is the reflection $t$ such that 
$t\gamma=\gamma s$.
That is, $t=\gamma s\gamma^{-1}$.
The map $\delta$, mapping each element to its left descent set, is a lattice homomorphism~\cite{Barbut} from the weak order on $W$ to the 
lattice of subsets of $S$ (see also~\cite[Section 9]{congruence}).
The maps $x\mapsto w_0x$ and $x\mapsto xw_0$ are anti-automorphisms, so that $x\mapsto w_0xw_0$ is an automorphism.

Given a subset $K\subseteq S$, let $W_K$ denote the subgroup of $W$ generated by $K$.
Then $(W_K,K)$ is a Coxeter system and $W_K$ is called a {\em parabolic subgroup} of $W.$\
In many sources, $W_K$ is called a {\em standard parabolic subgroup} and any conjugate of $W_K$ is called a parabolic subgroup, but in this paper, all parabolic subgroups are assumed to be standard.
The elements of $W_K$ form a lower interval in the weak order on $W,$\ and this interval in $W$ coincides with the weak order on $W_K$.
The maximal element of this interval is denoted $(w_0)_K$.

A finite group $W$ is a Coxeter group if and only if it has a representation as a group generated by orthogonal reflections in $\reals^n$ for some 
$n$.
In this case the reflections of $W$ are exactly the elements of $W$ which act as orthogonal reflections.
The set of reflecting hyperplanes for the reflections of $W$ is called the {\em Coxeter arrangement} $\A$ associated to $W.$\
The {\em regions} of $\A$ are the closures of the connected components of $\reals^n-\cup\A$.
A Coxeter arrangement is simplicial, in the sense that for each region $R$ of $\A$, choosing one normal vector to each facet of $R$ yields a 
linearly independent set of vectors. 
The regions of $\A$ are in one-to-one correspondence with the elements of $W,$\ and the identity element corresponds to one of the two 
antipodal regions whose facet hyperplanes are exactly the reflecting hyperplanes for the simple reflections.
The regions of $\A$, together with their faces, form a complete simplicial fan, which we denote $\F$.

Aside from a few general results, this paper primarily concerns the groups $A_n$ and $B_n$.
We now review the combinatorial interpretations of these groups in terms of permutations and signed permutations.

The Coxeter group $A_{n-1}$ is the symmetric group $S_n$ of permutations of $[n]:=\set{1,2,\ldots n}$.
The simple reflections are the adjacent transpositions $s_i:=(i,i+1)$ for $i\in[n-1]$, and the reflections are the transpositions.
The pairwise orders are $m(s_i,s_{i+1})=3$ for every $i\in[n-2]$ and $m(s_i,s_j)=2$ whenever $|i-j|>1$.
Thus the Coxeter diagram is a path with $n-1$ vertices and $n-2$ unlabeled edges.
We represent a permutation $x\in S_n$ in {\em one-line notation} $x_1x_2\cdots x_n$, meaning that $x$ maps~$i$ to $x_i$.
The left inversion set of $x$ is the set $\set{(x_j,x_i):1\le i<j\le n,x_i>x_j}$, and the left descent set is $\set{(x_j,x_i):1\le i<j\le n,x_i-1=x_j}$.
Cover relations in the weak order are transpositions of adjacent entries $x_i$ and $x_{i+1}$, and going up means switching $x_i$ and $x_{i+1}$
so as to put them out of numerical order.
The right descents of $x$ are the adjacent transpositions $(i,i+1)$ such that $x_i>x_{i+1}$.
The join-irreducibles of the weak order correspond to certain subsets $A\subseteq [n]$.
Defining $A^c:=[n]-A$,  $m=\min A$ and $M=\max A^c$, the join-irreducibles are the subsets with $M>m$.
Given such a subset, we construct the one-line notation for a join-irreducible permutation $\gamma$ by placing the elements of $A^c$ in 
ascending order followed by the elements of $A$ in ascending order.
The left reflection associated to $\gamma$ is $(m,M)$.

The maximal element $w_0$ of $S_n$ is the permutation $n(n-1)\cdots 1$.
We have $w_0x=(n+1-x_1)(n+1-x_2)\cdots(n+1-x_n)$ and $xw_0=x_nx_{n-1}\cdots x_1$.
Given permutations $x=x_1x_2\cdots x_m\in S_m$ and $y=y_1y_2\cdots y_n\in S_n$, say that the {\em pattern} $x$ 
{\em occurs in} $y$ if there are integers $1\le i_1<i_2<\cdots<i_m\le n$ such that  for all $1\le j<k\le m$ we have $x_j<x_k$ 
if and only if $y_{i_j}<y_{i_k}$.
Otherwise, say that $y$ {\em avoids} $x$.
For more information on patterns in permutations, see \cite{Wilf}.

The Coxeter arrangement corresponding to $S_n$ is most easily constructed in $\reals^n$, taking the hyperplanes normal to the vectors 
$e_i-e_j$ for $1\le j<i\le n$, where $e_i$ is the $i$th standard basis vector.
The identity permutation $12\cdots n$ corresponds to the region consisting of points $x$ with $x_1\le x_2\le\cdots\le x_n$, and the reflection 
$(j,i)$ fixes the hyperplane normal to $e_i-e_j$.

The Coxeter group $B_n$ is the group of {\em signed permutations}, meaning permutations $x$ of $\pm[n]:=\set{\pm1,\pm2,\ldots,\pm n}$ 
with the property that $x_{-i}=-x_i$ for every $i\in[n]$.
The simple reflections are $s_0:=(-1,1)$ and $s_i:=(i,i+1)(-i-1,-i)$ for $i\in[n-1]$.
The reflections are $(i,j)(-j,-i)$ for $j\neq-i$ and $(-i,i)$ for $i\in[n]$.
The pairwise orders are $m(s_i,s_{i+1})=3$ for every $i\in[n-1]$, $m(s_0,s_1)=4$ and $m(s_i,s_j)=2$ whenever $|i-j|>1$.
Thus the Coxeter diagram is a path with $n$ vertices and $n-1$ edges, each of which is unlabeled except for the first edge $s_0$---$s_1$, 
which is labeled $4$.
Signed permutations can be represented in {\em full (one-line) notation} $x_{-n},\ldots,x_{-1},x_1,\ldots,x_n$ where $x:i\mapsto x_i$,
or in {\em abbreviated notation} $x_1,\ldots,x_n$.
The left descent set of a signed permutation $x$ is $\set{(x_j,x_i)(-x_i,-x_j):i,j\in\pm[n],i<j,1\le x_i-1=x_j}$ together with $\set{(-1,1)}$
if the entry $1$ precedes the entry $-1$ in $x$.

Cover relations in the weak order can be described in terms of full notations for signed permutations as follows:
One type of cover transposes two symmetric pairs $x_i,x_{i+1}$ and $x_{-i-1},x_{-i}$ of adjacent entries in $x$.
The other type of cover transposes $x_{-1}$ and $x_1$.
In either case, going up in weak order means putting the entries out of numerical order.
The right descents of $x$ are the products $(i,i+1)(-i-1,-i)$ of adjacent transpositions for each $i\in[n-1]$ with $x_i>x_{i+1}$, and the 
transposition $(-1,1)$ whenever $x_1<0$.
The maximal element $w_0$ is $(-1)(-2)\cdots(-n)$.
The anti-automorphisms $x\mapsto w_0x$ and $x\mapsto xw_0$ both have the effect of negating each entry of a signed permutation $x$, so the 
automorphism $x\mapsto w_0xw_0$ is the identity.

The join-irreducibles of $B_n$ correspond to certain {\em signed subsets} $A$ of $[n]$.
A signed subset $A$ of $[n]$ is a subset of $\pm[n]$ such that for every $i\in[n]$, we have $\set{-i,i}\not\subseteq A$.
Given a signed subset $A$, let $-A:=\set{a:-a\in A}$, let $\pm A:=A\cup-A$, let the superscript ``$c$'' mean complementation in $\pm[n]$ 
and set $m:=\min A$.
Notice that the expression $-A^c$ is unambiguous, since $(-A)^c=-(A^c)$.
If $|A|=n$, set $M:=-m$ and otherwise set $M:=\max (\pm A)^c$.
Join-irreducible elements of the weak order on $B_n$ correspond to signed subsets $A$ with $M>m$.
Given a signed subset $A$ of $[n]$ with $M>m$, form the full notation for a join-irreducible signed permutation $\gamma$ by placing the 
elements of $-A$ in ascending order, followed by $(\pm A)^c$ in ascending order, then the elements of $A$ in ascending order.
The left reflection associated to $\gamma$ is $(m,M)(-M,-m)$ if $M\neq -m$ or $(m,-m)$ if $M=-m$.

There are at least two senses in which one can consider ``pattern avoidance'' in signed permutations, and both come into play in this paper.
First one might think of the full notation for a signed permutation in $B_n$ as an ordinary permutation of the set $\pm[n]$ and use the definition
of pattern avoidance for permutations given above.
We use phrases such as ``the full notation avoids the pattern\ldots'' to describe this kind of pattern avoidance, which appears with an 
additional variation in the description of the elements of the Cambrian lattices of type~B (Theorem~\ref{per B}).
Second, one can consider {\em signed pattern avoidance}, as defined in~\cite{signed}.
To any sequence $(a_1,a_2,\ldots,a_p)$ of nonzero integers with distinct absolute values, we associate a {\em standard signed permutation} 
$\st(a_1,a_2,\ldots,a_p)$.
This is the signed permutation $x\in B_p$ such that $x_i<0$ if and only if $a_i<0$ and $|x_i|<|x_j|$ if and only if $|a_i|<|a_j|$.
So for example $\st(7$-3-5$1)=4$-2-31.
If $x$ is a signed permutation, we say that a signed permutation $y$ {\em contains the signed pattern} $x$ if there is a subsequence of the entries
of $y$ whose standard signed permutation is $x$.
Otherwise, say that $y$ {\em avoids} $x$.

\section{Lattice congruences of the weak order}
\label{cong}
In this section we give background information about lattice congruences and prove Theorem~\ref{iso}, which classifies the isomorphisms and 
anti-isomorphisms among Cambrian lattices.
We then quote results of~\cite{congruence} which characterize the lattice of congruences of the weak order on a Coxeter group of type~A or~B\@.

Let $P$ be a finite poset with an equivalence relation $\Theta$ defined on the elements of $P$.
Given $a\in P$, let $[a]_\Theta$ denote the equivalence class of $a$.
The equivalence relation is an {\em order congruence} if:
\begin{enumerate}
\item[(i) ] Every equivalence class is an interval.
\item[(ii) ] The projection $\pidown:P\rightarrow P$, mapping each element $a$ of $P$ to the minimal element in $[a]_\Theta$, is 
order-preserving.
\item[(iii) ] The projection $\piup:P\rightarrow P$, mapping each element $a$ of $P$ to the maximal element in $[a]_\Theta$, is 
order-preserving.
\end{enumerate}
Define a partial order on the congruence classes by $[a]_\Theta\le[b]_\Theta$ if and only if there exists $x\in[a]_\Theta$ and $y\in[b]_\Theta$ 
such that $x\le_Py$.
The set of equivalence classes under this partial order is $P/\Theta$, the {\em quotient} of $P$ with respect to $\Theta$, which is isomorphic to
the induced subposet $\pidown(P)$.
The map $\piup$ maps $\pidown(P)$ isomorphically onto $\piup(P)$ with inverse $\pidown$.
For more information on order congruences and quotients, see~\cite{Cha-Sn,Order}.

Let $L$ be a finite lattice, with join and meet operations denoted by $\join$ and $\meet$ respectively.
Recall that an element $\gamma$ of $L$ is join-irreducible if and only if it covers exactly one element $\gamma_*$ of~$L$.
Meet-irreducible elements are defined dually.
Denote the set of join-irreducibles of~$L$ by $\Irr(L)$.
We also use $\Irr(L)$ to represent the induced subposet of~$L$ consisting of join-irreducible elements.  
A {\em lattice congruence} on a lattice~$L$ is an equivalence relation on the elements of~$L$ which respects joins and meets in the sense that
if $a_1\equiv a_2$ and $b_1\equiv b_2$ then $a_1\join b_1\equiv a_2\join b_2$ and similarly for meets.
Given a congruence $\Theta$ on $L$, there is a well-defined meet and join on equivalence classes, namely, 
$[a]_\Theta\join[b]_\Theta=[a\join b]_\Theta$, and similarly for meets.
The resulting lattice of equivalence classes is the quotient $L/\Theta$.
An equivalence relation $\Theta$ on a lattice~$L$ is an order congruence if and only if it is a lattice congruence, and the two definitions of the 
quotient of~$L$ mod $\Theta$ coincide.

For a cover relation $x\covered y$, if $x\equiv y\mod\Theta$ then we say $\Theta$ {\em contracts} the edge $x\covered y$.
For an element $y$, if there exists an edge $x\covered y$ contracted by $\Theta$, we say $\Theta$ contracts $y$.
In particular $\Theta$ contracts a join-irreducible $\gamma$ if and only if $\gamma\equiv \gamma_*$.

Let $\Con(L)$ be the set of congruences of~$L$ partially ordered by refinement, which is a distributive lattice~\cite{Fun-Nak}, and thus is 
uniquely determined by the subposet $\Irr(\Con(L))$.
The join-irreducible lattice congruences (elements of $\Irr(\Con(L))$) can be characterized in terms of the join-irreducible elements of $L$ as follows.
The meet in $\Con(L)$ is intersection of congruences (considering the congruences as binary relations).
Given a join-irreducible $\gamma$ of~$L$, let $\Cg(\gamma)$ be the meet of all congruences contracting $\gamma$.
Since the meet is intersection, $\Cg(\gamma)$ also contracts $\gamma$.
Thus $\Cg(\gamma)$ is the unique smallest lattice congruence contracting $\gamma$.
Using a characterization of the join in $\Con(L)$ found in~\cite{Gratzer}, one sees that $\Cg(\gamma)$ is join-irreducible as an element of $\Con(L)$.
In fact, all of the join-irreducible elements of $\Con(L)$ are of the form $\Cg(\gamma)$ for $\gamma$ join-irreducible in $L$.
(In other words, a lattice congruence is determined by the set of join-irreducibles $\gamma$ it contracts~\cite[Section II.3]{Free Lattices}.)
Although the map $\Cg:\Irr(L)\to\Irr(\Con(L))$ is onto, it need not be one-to-one.
A lattice~$L$ is {\em congruence uniform} if $\Cg$ is a bijection and if a dual statement about meet-irreducibles holds as 
well~\cite{cong norm}.
In particular, when $L$ is a congruence uniform lattice, $\Irr(\Con(L))$ can be thought of as a partial order on the join-irreducibles via the bijection $\Cg$.

A lattice homomorphism is a map of lattices which respects the join and meet operations in the usual algebraic sense.
The fibers of a lattice homomorphism form a lattice congruence, and given a congruence $\Theta$, the canonical map $L\to L/\Theta$ is a 
homomorphism.
If $\eta:L\to L_1$ and $\zeta:L\to L_2$ are lattice homomorphisms, we say $\zeta$ {\em factors through} $\eta$ if there is a lattice
homomorphism $\zeta':L_1\to L_2$ such that $\zeta=\zeta'\circ\eta$.
Let $\Theta$ and $\Phi$ be the lattice congruences associated to $\eta$ and $\zeta$.
If $\Theta\le\Phi$ in $\Con(L)$ then $\zeta$ factors through $\eta$.
In this case, $\Phi$ can be thought of as a congruence on $L/\Theta$.
Thus for any congruence $\Theta$ on $L$, $\Irr(\Con(L/\Theta))$ is the filter in $\Irr(\Con(L))$ consisting of join-irreducibles of~$L$ not 
contracted by $\Theta$.

We omit the easy proofs of the following facts.
\begin{prop}
\label{auto sublattice}
Let $L$ be a finite lattice and $G$ a group acting on $L$ by automorphisms.
Then $L^G$, the subposet of~$L$ consisting of elements fixed by every $g\in G$, is a sublattice of~$L$.
\end{prop}

\begin{prop}
\label{sub congruence}
Let $L$ be a finite lattice, $L'$ a sublattice and $\Theta$ a congruence on~$L$.
Then the restriction of $\Theta$ to $L'$ is a congruence on $L'$.
\end{prop}

For the remainder of the section, we focus specifically on congruences of the weak order on a finite Coxeter group $W.$\
The weak order on $W$ is a congruence uniform lattice~\cite{boundedref} so $\Irr(\Con(W))$ is a partial order on the join-irreducibles of $W.$\
Since a parabolic subgroup $W_K$ is a lower interval in the weak order on $W,$\ join-irreducibles in $W_K$ are in particular join-irreducibles 
of $W.$\
The following is~\cite[Lemma 6.12]{congruence}, specialized to the weak order on a Coxeter group.
\begin{lemma}
\label{par lemma}
If $W_K$ is a parabolic subgroup of $W,$\ then $\Irr(\Con(W_K))$ is an order filter in $\Irr(\Con(W))$.
\end{lemma}
That is, if $\gamma\in\Irr(\Con(W_K))\subseteq\Irr(\Con(W))$ and $\gamma'\ge \gamma$, then we have $\gamma'\in\Irr(\Con(W_K))$.
The {\em degree} $\deg(\gamma)$ of a join-irreducible $\gamma$ is $|K|$ for the smallest $K$ such that $\gamma\in W_K$.
By Lemma~\ref{par lemma}, if $\gamma\le \gamma'$ in $\Irr(\Con(W))$ then $\deg(\gamma)\ge \deg(\gamma')$.
A congruence $\Theta$ is {\em homogeneous of degree $k$} if there is some set $C$ of join-irreducibles of degree $k$ such that $\Theta$
is the smallest congruence contracting all join-irreducibles in $C$.

Let $\alpha$ denote the map $w\mapsto ww_0$ on $W.$\
If $\Theta$ is a congruence on the weak order on $W,$\ let $\alpha(\Theta)$ be the antipodal congruence, defined by 
$\alpha(x)\equiv \alpha(y)\mod\alpha(\Theta)$ if and only if $x\equiv y\mod\Theta$.
Then $\alpha$ induces an anti-isomorphism from $W/\Theta$ to $W/(\alpha(\Theta))$.
The following is part of ~\cite[Proposition 6.13]{congruence}.
\begin{proposition}
\label{dual cong}
Let $\gamma=sts\cdots$ be a reduced word for a degree-two join-irreducible in a Coxeter group $W.$\
Then $\gamma':=\gamma_*\cdot(w_0)_{\set{s,t}}$ is a degree-two join-irreducible with reduced word of the form $tst\cdots$ with 
$l(\gamma')=m(s,t)-l(\gamma)+1$ and $\gamma$ is contracted by $\Theta$ if and only if $\gamma'$ is contracted by $\alpha(\Theta)$.
\end{proposition}

We apply these facts to the problem of determining when two Cambrian lattices are isomorphic or anti-isomorphic. 
An isomorphism of directed diagrams is a bijection on the vertices which preserves labels and directions of edges, while an anti-isomorphism 
preserves edge-labels and reverses arrows.
\begin{theorem}
\label{iso}
An isomorphism (respectively anti-isomorphism) of directed diagrams induces an isomorphism (respectively anti-isomorphism) of Cambrian
lattices.
An isomorphism (respectively anti-isomorphism) of Cambrian lattices restricts to an isomorphism (respectively anti-isomorphism) of directed
diagrams.
\end{theorem}
\begin{proof}
The atoms of the weak order on $W$ are join-irreducibles of degree 1.
The Cambrian congruence $\Theta(\G)$ is the homogeneous congruence of degree 2 generated by contracting the  join-irreducibles with 
reduced words of the form $tst\cdots$ with lengths weakly between $2$ and $m(s,t)-1$, for every edge $s\to t$ of $\G$.
It is easily checked that $\Irr(\Con(W_{\set{s,t}}))$ is the partially ordered set in which $\set{s,t}$ is an antichain, the remaining 
join-irreducibles form an antichain, and $s>\gamma$ and $t>\gamma$ for every join-irreducible $\gamma\not\in\set{s,t}$.
Since, by Lemma~\ref{par lemma}, $\Irr(\Con(W_{\set{s,t}}))$ is an order filter in $\Irr(\Con(W))$, the congruence $\Theta(\G)$ does not 
contract any other join-irreducibles in $W_{\set{s,t}}$ besides those listed above, and furthermore does not contract any other elements of 
$W_{\set{s,t}}$ besides those join-irreducibles.
The maximal element of $W_{\set{s,t}}$ is the join of $s$ and $t$ in $W,$\ and since it is not contracted, it is also the join of $s$ and $t$ in 
$\C(\G)$.
Therefore, given an abstract partial order known to be isomorphic to some Cambrian lattice $\C(\G)$, we can reconstruct $\G$ as follows:
The vertices of $\G$ are the atoms of the lattice.
For each pair $\set{s,t}$ of distinct atoms, consider the interval $[1,s\join t]$.
If this interval has 4 elements, then $s$ and $t$ do not form an edge in $G$.
Otherwise the interval consists of two chains from $1$ to $s\join t$, one containing 3 elements, and the other containing $m+1$ elements for 
some $m>2$.
Without loss of generality, let $s$ be contained in the larger chain.
Then $\G$ has an edge labeled $m$ pointing from $s$ to $t$.

We have established the second statement for isomorphisms, and the first statement for isomorphisms is immediate from the definition.
The statements for anti-isomorphisms follow by Proposition~\ref{dual cong}.
\end{proof}

The proof of Theorem~\ref{iso} also shows that for any $K\subseteq S$, the restriction to $W_K$ of the Cambrian congruence 
$\Theta(\G)$ is $\Theta(\G_K)$, where $\G_K$ is the directed subgraph of $\G$ induced by $K$.
In particular, the lower intervals in the Cambrian lattice $\C(\G)$ arising from (standard) parabolic subgroups of $W$ are Cambrian lattices.

In Section 9 of~\cite{congruence} it is observed that the congruence associated to the left descent map on $W$ contracts every 
join-irreducible of degree $>1$.
Thus the descent map factors through any lattice homomorphism $\eta$ whose associated congruence contracts no join-irreducibles of degree 1.
In particular, the descent map factors through the canonical homomorphism from $W$ to the Cambrian lattice $\C(\G)$.

In~\cite{congruence}, the author determined the congruence lattice $\Con(W)$ of the weak order on $W$ for $W=S_n$ or $W=B_n$.
Since the weak order is congruence uniform, $\Irr(\Con(W))$ is described as a partial order on the join-irreducibles of $W.$\
Recall from Section~\ref{weak} that join-irreducibles of $S_n$ correspond to certain subsets $A\subseteq [n]$, and recall the associated notation
including $m:=\min A$ and $M:=\max A^c$.
The degree of a join-irreducible $\gamma$ is $M-m$.
We state~\cite[Theorem 8.1]{congruence} as follows, where the transitive closure of a directed graph is interpreted as a poset via the 
convention that $x\to y$ corresponds to $x\ge y$.
\begin{theorem}
\label{A shard}
The poset $\Irr(\Con(S_n))$ is the transitive closure of the directed graph in which $\gamma_1\to \gamma_2$ if and only if the corresponding subsets $A_1$ 
and $A_2$ satisfy one of the following:
\begin{enumerate}
\item[(i) ]$A_1\cap[1,M_1)=A_2\cap[1,M_1)$ and $M_2>M_1$, or
\item[(ii) ]$A_1\cap(m_1,n]=A_2\cap(m_1,n]$ and $m_2<m_1$.
\end{enumerate}
\end{theorem}
Using Theorem~\ref{A shard} one can determine~\cite[Theorem 8.2]{congruence} the cover relations in $\Irr(\Con(S_n))$ and conclude that
$\Irr(\Con(S_n))$ is dually ranked by degree, in the sense that if $\gamma_1\covers\gamma_2$ then $\deg(\gamma_1)=\deg(\gamma_2)-1$.

Recall from Section~\ref{weak} that join-irreducibles of $B_n$ correspond to certain signed subsets $A$ of $ [n]$, and recall the 
definitions of $m$ and $M$ for signed subsets.
The following lemma describes the lower-degree join-irreducibles of $B_n$.
\begin{lemma}
\label{B degrees}
For a join-irreducible $\gamma$ of $B_n$ with associated signed subset $A$:
\begin{enumerate}
\item[(i) ] $\deg(\gamma)=1$ if and only if $|(m,M)\cap(\pm[n])|=0$, and
\item[(ii) ] $\deg(\gamma)=2$ if and only if either $|(m,M)\cap(\pm[n])|=1$ or $(m,M)=(-2,2)$.
\end{enumerate}
\end{lemma}
\begin{proof}
Since $M>0$, we have $|(m,M)\cap(\pm[n])|=0$ if and only if $(m,M)=(-1,1)$ or $(m,M)=(i,i+1)$ for some $i\in[n-1]$.
Each of these choices of $m$ and $M$ yields a unique signed subset corresponding to an element of degree one.
We have $|(m,M)\cap(\pm[n])|=1$ if and only if $(m,M)=(-2,1)$, $(-1,2)$ or $(i,i+2)$ for $i\in[n-2]$.
It is easily checked that each of the first two pairs $(m,M)$ corresponds to a single join-irreducible in the parabolic subgroup 
$W_{\set{s_0,s_1}}$, and that for $(m,M)=(i,i+2)$, there are two corresponding join-irreducibles, both in $W_{\set{s_i,s_{i+1}}}$.
Finally, for $(m,M)=(-2,2)$ there are two join-irreducibles, both in $W_{\set{s_0,s_1}}$.
We have accounted for all of the join-irreducibles of degree 1 or 2.
\end{proof}

The characterization of $\Irr(\Con(B_n))$ is more complicated than the \mbox{type-A} analog, but still serviceable for our purposes.
Let $A_1$ and $A_2$ be signed subsets of $[n]$ with $M_1,M_2,m_1,m_2$ as defined above.
To simplify the statement of the theorem, we first list some conditions on $A_1$ and $A_2$.
The conditions are labeled by letters q, f and r, and these letters had meanings in~\cite{congruence} which we need not explain here.
We keep the labelings for the sake of consistency between the two papers. 
First, conditions (q1) through (q6):

\begin{tabular}{ll}
(q1)&$-m_1=M_1<M_2=-m_2$.\\
(q2)&$-m_2=M_2=M_1>m_1>0$.\\
(q3)&$M_2=M_1>m_1>m_2\neq -M_2$.\\
(q4)&$M_2>M_1>m_1=m_2\neq -M_2$.\\
(q5)&$-m_2=M_1>m_1>-M_2\neq m_2$.\\
(q6)&$-m_2>M_1>m_1=-M_2\neq m_2$.
\end{tabular}

\noindent
Next we have condition (f) which is a combination of three conditions, and which depends on a parameter $a$ in $\pm[n]$.
Say $A_2$ satisfies condition (f\,:\,$a$) if one of the following holds:

\begin{tabular}{ll}
(f1\,:\,$a$)&$a\in A_2$.\\
(f2\,:\,$a$)&$a\in A_2^c-\set{-M_2,-m_2}$ and $-a\not\in A_2\cap(m_2,M_2)$.\\
(f3\,:\,$a$)&$a\in\set{-M_2,-m_2}$ and $(\pm A_2)^c\cap(m_2,M_2)\cap(-M_2,-m_2)=\emptyset$.
\end{tabular}

\noindent
Finally, we have conditions (r1) and (r2):

\begin{tabular}{ll}
(r1)&$A_1\cap(m_1,M_1)=A_2\cap(m_1,M_1)$.\\
(r2)&$A_1\cap(m_1,M_1)=-A_2^c\cap(m_1,M_1)$.
\end{tabular}

\begin{theorem}\cite[Theorem 7.3]{congruence}
\label{B shard}
The poset $\Irr(\Con(B_n))$ is the transitive closure of the directed graph in which $\gamma_1\to \gamma_2$ if and only if the corresponding subsets $A_1$ 
and $A_2$ satisfy one or more of the following combinations of conditions:
\begin{enumerate}
\item[1. ] {\rm(q1)} and {\rm(r1)}.
\item[2. ] {\rm(q2)} and {\rm(r1)}.
\item[3. ] {\rm(q3)}, {\rm(f\,:\,$m_1$)} and {\rm(r1)}.
\item[4. ] {\rm(q4)}, {\rm(f\,:\,$M_1$)} and {\rm(r1)}.
\item[5. ] {\rm(q5)}, {\rm(f\,:\,$\,-m_1$)} and {\rm(r2)}.
\item[6. ] {\rm(q6)}, {\rm(f\,:\,$\,-M_1$)} and {\rm(r2)}.
\end{enumerate}
\end{theorem}

\section{The Type-A Associahedron}
\label{A assoc}
In this section, we review an iterated fiber polytope construction due to Billera and Sturmfels~\cite{Fiber,Iterated} which yields a map from 
the permutohedron to the associahedron.
For details on fiber polytopes, see~\cite{Fiber,Iterated}.

Given polytopes $P$ and $Q$ and a linear surjection $P\mapname{\varphi}Q$, the fibers of $\varphi$ form a {\em bundle of polytopes} over $Q$,
in a sense made precise in~\cite{Fiber}.
Billera and Sturmfels define a {\em Minkowski integral} over polytope bundles, and the Minkowski integral, subject to some normalization, of
the bundle of fibers of $\varphi$ is called the {\em fiber polytope} $\Sigma(P\mapname{\varphi}Q)$.
Given a tower $P\mapname{\varphi}Q\mapname{\rho}R$ of surjective linear maps of polytopes, the fiber polytope 
$\Sigma(Q\mapname{\rho}R)$ is the image of $\Sigma(P\mapname{\rho\circ\varphi}R)$ under the map $\varphi$.
The {\em iterated fiber polytope} $\Sigma(P\mapname{\varphi}Q\mapname{\rho}R)$ is defined to be the fiber polytope 
$\Sigma\left(\Sigma(P\mapname{\rho\circ\varphi}R)\mapname{\varphi}\Sigma(Q\mapname{\rho}R)\right)$.
Theorem 2.1 of~\cite{Iterated} states that the normal fan of $\Sigma(P\mapname{\varphi}Q\mapname{\rho}R)$ refines that of 
$\Sigma(P\mapname{\varphi}Q)$.

Consider the tower of surjective linear maps of polytopes 
\[\Delta^{n+1}\mapname{\varphi} Q_{n+2}\mapname{\rho} I,\]
where $I$ is a 1-dimensional polytope, $Q_{n+2}$ is a polygon with $n+2$ vertices, and 
$\Delta^{n+1}$ is the $(n+1)$-dimensional simplex whose vertices are the coordinate vectors $ e_0, e_1,\ldots, e_{n+1}$ 
in $\reals^{n+2}$.
We sometimes refer to these polytopes  simply as $\Delta$ and $Q$.
Let $v_i:=\varphi( e_i)$ and $a_i:=\rho(v_i)$ and suppose that $a_0<a_1<\cdots<a_{n+1}$.
Let $f$ be a non-trivial linear functional on $\ker\rho$.
We may as well take $\rho$ to be an orthogonal projection of $Q$ onto the line segment whose endpoints are $v_0$ and $v_{n+1}$ and think of 
$f$ as giving the positive or negative ``height'' of each vertex of $Q_{n+2}$ above that line segment.
We abbreviate $f_i:=f(v_i)$ and use the shorthand $\upwide{i}$ to denote an $i\in[n]$ with $f_i\ge 0$.
In this case we call $v_i$ an {\em up vertex} and~$i$ an {\em up index}.
Similarly, $\downwide{i}$ denotes an $i\in[n]$ with $f_i\le 0$, called a {\em down index}.
Thus for example the phrase ``Let  $\downwide{i}\in[n]$'' means ``Let $i\in[n]$ have $f_i\le 0$.''
For any $H\subseteq[n]$, let $\down{H}=\set{\downwide{i}\in H}$ and let $\up{H}=\set{\upwide{i}\in H}$.

The fiber polytope $\Sigma(\Delta\mapname{\rho\circ\varphi}I)$ is a cube whose vertices correspond to triangulations of the point 
configuration $\set{a_0,a_1,\ldots,a_{n+1}}$.
Such a triangulation can be thought of as a subset $H=\set{i_1,i_2,\ldots,i_k}\subseteq [n]$ where $H$ is the set of points, other than the 
endpoints of $I$, appearing as vertices of the triangulation.

The iterated fiber polytope $\Sigma(\Delta\to Q\to I)$ is a polytope of whose vertices correspond to $f$-monotone paths in 
$\Sigma(\Delta\mapname{\rho\circ\varphi}I)$.
These paths are permutations, and $\Sigma(\Delta\to Q\to I)$ is combinatorially isomorphic to the $(S_n)$-permutohedron.
Specifically, an $f$-monotone path $\down{[n]}=H_0\to H_1\to\cdots\to H_n=\up{[n]}$ is associated to the permutation 
$x_1x_2\cdots x_n$ where $x_i$ is the unique element in the symmetric difference of $H_i$ and $H_{i-1}$.
Two such monotone paths are connected by an edge in $\Sigma(\Delta\to Q\to I)$ if they differ in only one vertex.
Thus edges correspond to cover relations in the weak order.
To simplify the statements of some results, it is sometimes convenient to think of $x$ as the word 
$x_0x_1x_2\cdots x_nx_{n+1}$ where $x_0$ is always 0, and $x_{n+1}$ is always $n+1$.
Furthermore, it is useful to remember that 0 is both a down vertex $\down{0}$ and an up vertex $\up{0}$ and that $n+1$ is both 
$\down{n+1}$ and $\up{n+1}$.
Notice also that reversing the sign of $f$ corresponds to reversing the direction of each $f$-monotone path, which corresponds to reversing 
the order of the entries of each permutation.
Thus reversing the sign of $f$ is the anti-automorphism $x\mapsto xw_0$ of weak order.
Re-indexing the points $v_i$ and $a_i$ by $i\mapsto n-i+1$ is the anti-automorphism $x\mapsto w_0x$.
These anti-automorphisms and their composition $x\mapsto w_0xw_0$ simplify the proofs of some results.

The fiber polytope $\Sigma(\Delta^{n+1}\mapname{\varphi} Q_{n+2})$ is combinatorially isomorphic to the $S_n$-associahedron,
whose vertices are the triangulations of $Q$.
We think of a triangulation as a maximal set of non-intersecting diagonals of $Q$.
Faces of the associahedron correspond to sets of partial triangulations, or in other words triangulations from which some number of diagonals 
have been removed.
Edges are partial triangulations that are missing one diagonal.
Such a partial triangulation defines a unique quadrilateral, and can be completed to a triangulation by inserting either of 
the two diagonals of the quadrilateral.
A {\em diagonal flip} is the process of removing an edge from a triangulation and completing the resulting partial triangulation to a 
triangulation by inserting the other diagonal of the quadrilateral.
As mentioned above in the general case, the normal fan of $\Sigma(\Delta\to Q\to I)$ refines that of $\Sigma(\Delta\to Q)$.
In particular, there is a map $\eta$ sending a maximal normal cone $C$ of $\Sigma(\Delta\to Q\to I)$ to the maximal normal cone of 
$\Sigma(\Delta\to Q)$ containing $C$.

\begin{figure}[ht]
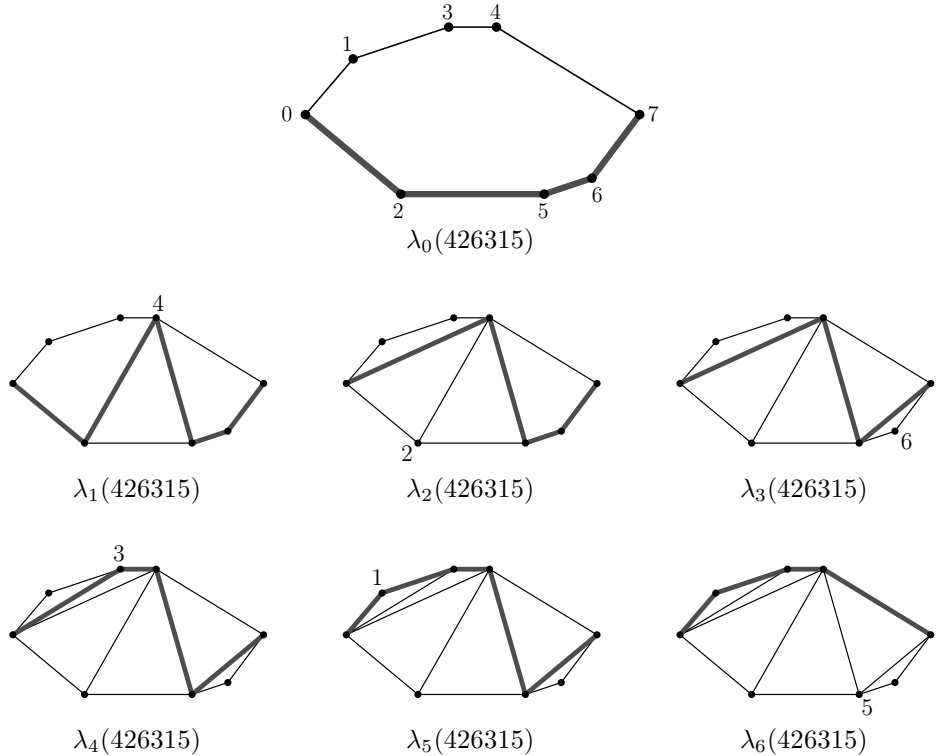

\centerline{$\begin{array}{c}
\scalebox{0.6}{\epsfbox{lambda0.ps}}\\
\lambda_0(426315)\\[15 pt]
\end{array}
$}   
\centerline{$\begin{array}{ccc}
\begin{array}{c}\scalebox{0.45}{\epsfbox{lambda1.ps}}\\[2 pt]\lambda_1(426315)\end{array}  &
\begin{array}{c}\scalebox{0.45}{\epsfbox{lambda2.ps}}\\[2 pt]\lambda_2(426315)\end{array}  &
\begin{array}{c}\scalebox{0.45}{\epsfbox{lambda3.ps}}\\[2 pt]\lambda_3(426315)\end{array}  \\[50 pt]
\begin{array}{c}\scalebox{0.45}{\epsfbox{lambda4.ps}}\\[2 pt]\lambda_4(426315)\end{array}  &
\begin{array}{c}\scalebox{0.45}{\epsfbox{lambda5.ps}}\\[2 pt]\lambda_5(426315)\end{array}  &
\begin{array}{c}\scalebox{0.45}{\epsfbox{lambda6.ps}}\\[2 pt]\lambda_6(426315)\end{array}  
\end{array}$}
\caption{The map $\eta$.}
\label{eta}
\end{figure}

The map $\eta$ takes a permutation $x=x_1x_2\cdots x_n$ to a triangulation of $Q$, and has a characterization in terms of polygonal 
paths.
The diagonals of the triangulation arise as a union of polygonal paths $\lambda_0(x),\lambda_1(x),\ldots,\lambda_n(x)$ in $Q$ such that each 
vertex of each path is a vertex of $Q$, and such that each path visits vertices in the order given by their subscripts.
Specifically, if $x$ is the permutation associated to the monotone path $H_0\to H_1\to\cdots\to H_n$, then $\lambda_i(x)$ visits the vertices 
$\set{v_j:j\in H_i}$ in the order given by their subscripts.
Alternately, one can let $\lambda_0(x)$ be the path from $v_0$ to $v_{n+1}$ passing through the points $v_i$ for $\down{i}\in [n]$ and
define $\lambda_i(x)$ recursively:
If $x_i$ is $\down{x_i}$, define $\lambda_i(x)$ by deleting $v_{x_i}$ from the list of vertices visited by $\lambda_{i-1}(x)$, and 
if $x_i$ is $\up{x_i}$, define $\lambda_i(x)$ by adding $v_{x_i}$ to the list of vertices visited by $\lambda_{i-1}(x)$.
Figure~\ref{eta} illustrates, for the case $n=6$, the map $\eta$ applied to the permutation $x=426315$, where $Q$ is a polygon with $\down{[n]}=\set{2,5,6}$ and $\up{[n]}=\set{1,3,4}$.

The combinatorics of the map $\eta$ from permutations to triangulations of $Q$ derive from the sign of $f$ on each vertex of $Q$.
Thus to be more exact, we should name the map $\eta_f$.
Usually, however, the choice of $f$ is fixed, so we drop the subscript $f$ unless we want to emphasize the fact that $f$ can vary.

\section{Properties of the map $\eta$}
\label{eta A}
\begin{figure}[ht]
\centerline{\scalebox{.6}{\epsfbox{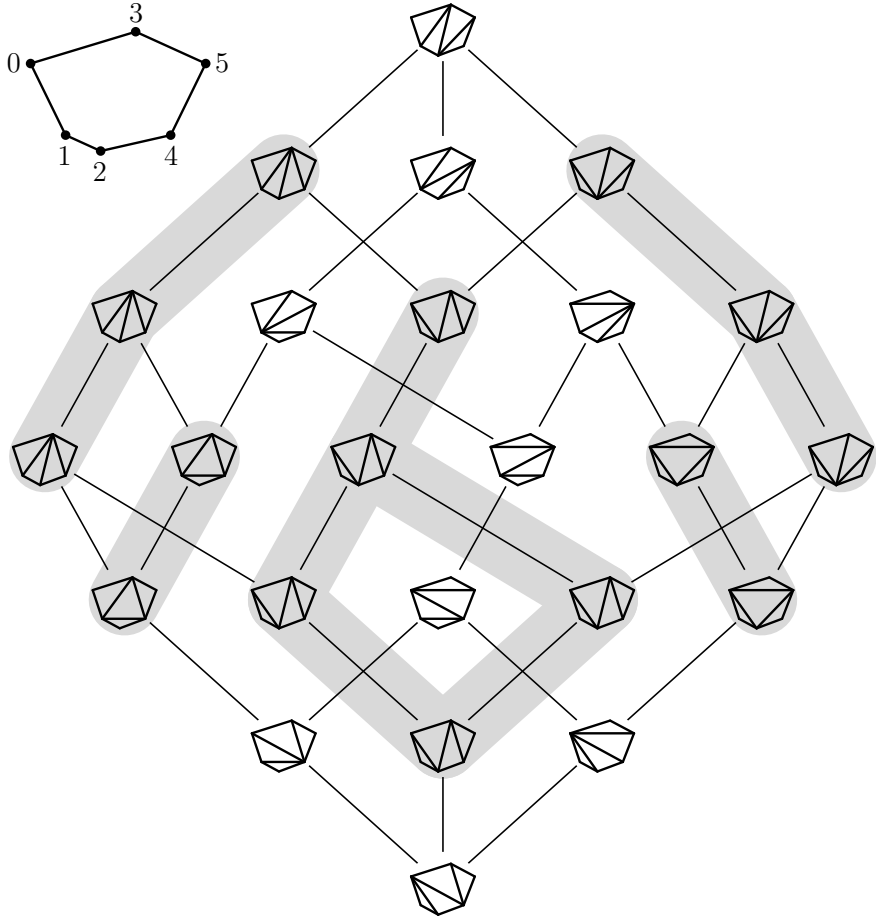}}}   
\caption{Fibers of $\eta$.}
\label{A3tri}
\end{figure}

In this section we to prove the following theorem, and in the process obtain some useful results about the map $\eta$.
\begin{theorem}
\label{A congruence}
The fibers of $\eta_f$ are the congruence classes of a lattice congruence $\Theta_f$ on the weak order on $S_n$.
\end{theorem}
Figure~\ref{A3tri} illustrates Theorem~\ref{A congruence} by showing the Hasse diagram of the weak order on $S_4$, as drawn in Figure~\ref{A3camb_cong}, with each element $x\in S_4$ replaced by $\eta_f(x)$.
Here $f$ defines up-indices $\up{[n]}=\set{3}$ and down-indices $\down{[n]}=\set{1,2,4}$.
The shaded edges are edges $x\covered y$ such that $\eta(x)=\eta(y)$.
These shaded edges define a lattice congruence which coincides with the Cambrian congruence pictured in Figure~\ref{A3camb_cong}
(see Theorem~\ref{cambrian A}).

The proof of Theorem~\ref{A congruence} occupies the remainder of this section.
To begin with, the fact that $\eta$ arises in connection with a refinement of the normal fan means that the inverse image under $\eta$ of a 
vertex of $\Sigma(\Delta\to Q)$ is a connected subgraph of the 1-skeleton of $\Sigma(\Delta\to Q\to I)$.
Since the 1-skeleton of the permutohedron is isomorphic to the Hasse diagram of the weak order, we have:
\begin{lemma}
\label{connected}
The fibers of $\eta$ are connected subgraphs of the Hasse diagram of the weak order on $S_n$.
\end{lemma}

We use a variant of pattern-avoidance to characterize the fibers of $\eta$.
Except in special cases, this is not pattern-avoidance in the sense of Section~\ref{weak}, because it requires keeping track of 
whether certain indices $i\in[n]$ are up or down indices.
Say that $x$ contains the pattern $\up{2}31$ if there exist $1\le i<j<k\le n$ with $x_k<\up{x_i}<x_j$.
This means $x_k<x_i<x_j$ and $f_{x_i}>0$.  
No conditions are placed on $f_{x_j}$ or $f_{x_k}$.
In the special case $k=j+1$, altering $x$ by switching the entries $x_j$ and $x_{j+1}$ is called a $\up{2}31\to\up{2}13$-move, and this move 
goes down by a cover in the weak order.
Similarly, $x$ contains the pattern $31\down{2}$ if there exist $1\le i<j<k\le n$ with $x_j<\down{x_k}<x_i$.
If $j=i+1$, then altering $x$ by switching the entries $x_i$ and $x_{i+1}$ is called a $31\down{2}\to 13\down{2}$-move, which is also a cover
in weak order.
One can similarly define the notion of permutations containing $\up{2}13$ and $13\down{2}$.

As mentioned in Section~\ref{A assoc}, the edges of the $S_n$-associahedron correspond to diagonal flips.
Given two triangulations $T_1$ and $T_2$ differing by a diagonal flip, we denote by $T_1-T_2$ the unique diagonal in $T_1$ which is not 
a diagonal in $T_2$.
The slope of a diagonal refers to the usual slope, relative to the convention that the positive horizontal direction is the direction of a ray from
$v_0$ through $v_{n+1}$ and the positive vertical direction is the positive direction of the functional $f$.

\begin{prop}
\label{cover moves}
Suppose $x$ and $y$ are permutations of $[n]$ and $x\covered y$ in weak order.
Then the following are equivalent:
\begin{enumerate}
\item[(i) ] $\eta(x)=\eta(y)$.
\item[(ii) ] $x$ is obtained from $y$ by a $\up{2}31\to\up{2}13$-move or a $31\down{2}\to 13\down{2}$-move.
\end{enumerate}
Furthermore, if $\eta(x)\neq\eta(y)$, then the two triangulations differ by a diagonal flip, and the slope of the diagonal $\eta(x)-\eta(y)$
is less than the slope of the diagonal $\eta(y)-\eta(x)$.
\end{prop}
\begin{proof}
Since $x\covered y$, we can write $y=x_1x_2\cdots x_{j-1}x_{j+1}x_jx_{j+2}x_{j+3}\cdots x_n$ for some $j\in[n-1]$ with $x_j<x_{j+1}$.
We have $\lambda_k(x)=\lambda_k(y)$ for every $k$ except $k=j$.

First, we  show that (ii) implies (i).
If $x$ is obtained from $y$ by a $\up{2}31\to\up{2}13$-move, let~$i$ be the position of the entry described by the $\up{2}$ in this move.
(The reader may want to consider Figure~\ref{eta}, taking $x=426315$, $y=462315$ and $i=1$.)
For every $k\ge i$, the path $\lambda_k(x)$ consists of a path from $v_0$ to $v_{x_i}$ composed with a path from $v_{x_i}$ to $v_{n+1}$.
Also for every $k>i$, if $x_k<x_i$ then the path $\lambda_k(x)$ agrees, from $v_{x_i}$ to $v_{n+1}$, with the path $\lambda_{k-1}(x)$, 
and if $x_k>x_i$ then the path $\lambda_k(x)$ agrees, from $v_0$ to $v_{x_i}$, with the path $\lambda_{k-1}(x)$.
In particular, $\lambda_j(x)$ agrees with $\lambda_{j+1}(x)$ from $v_0$ to $v_{x_i}$ and agrees with $\lambda_{j-1}(x)$ from 
$v_{x_i}$ to $v_{n+1}$.
Thus $\eta(x)$ can be obtained from the union of the paths $\lambda_k(x)$ for $k\neq j$.
Reasoning similarly, we see that $\eta(y)$ can be obtained from the union of the paths $\lambda_k(y)$ for $k\neq j$.
Since $\lambda_k(x)=\lambda_k(y)$ for $k\neq j$, we have $\eta(x)=\eta(y)$.
The case of a $31\down{2}\to 13\down{2}$-move is symmetrical by $x\mapsto w_0xw_0$, thus completing the proof that (ii) implies (i).

Suppose that $x$ is not obtained from $y$ by a $\up{2}31\to\up{2}13$-move or a $31\down{2}\to 13\down{2}$-move.
Thus every $\up{x_i}$ with $i<j$ has either $x_i<x_j$ or $x_i>x_{j+1}$, and every $\down{x_i}$ with $i>j+1$ has either 
$x_i<x_j$ or $x_i>x_{j+1}$.
In particular the path $\lambda_{j-1}(x)=\lambda_{j-1}(y)$ has no vertices $v_k$ with $x_j<k<x_{j+1}$.
The same is true of the paths $\lambda_j(x)$, $\lambda_j(y)$ and $\lambda_{j+1}(x)=\lambda_{j+1}(y)$.
Breaking into four cases, depending on whether each of $x_j$ and $x_{j+1}$ is an up index or a down index, one easily checks the following:
The union of the paths $\lambda_{j-1}(x)$, $\lambda_j(x)$ and $\lambda_{j+1}(x)$ encloses a quadrilateral cut in two by one of its diagonals.
The union of $\lambda_{j-1}(y)$, $\lambda_j(y)$ and $\lambda_{j+1}(y)$ encloses the same quadrilateral cut in two by the other diagonal.
In particular $\eta(x)\neq\eta(y)$, and the statement about slopes is also easily checked in each of the four cases.
\end{proof}

\begin{prop}
\label{convex}
For any permutation $x$, the fiber $\eta^{-1}(\eta(x))$ is order-convex in the weak order.
\end{prop}
\begin{proof}
For convenience, let $\slo(T)$ denote the sum of the slopes of the diagonals of a triangulation $T$.
We need to show that if $x<y<z$ and $\eta(x)=\eta(z)$ then we must have $\eta(y)=\eta(z)$.
If not, then by repeated applications of Proposition~\ref{cover moves} we find that the $\slo(\eta(x))<\slo(\eta(y))<\slo(\eta(z))$.
This contradicts the fact that $\eta(x)=\eta(z)$, thus proving that $\eta(y)=\eta(z)$.
\end{proof}

\begin{lemma}
\label{adjacent}
If $x\in S_n$ contains either the pattern $\up{2}31$ or $31\down{2}$, then $x$ contains some instance of that pattern such that the 
``3'' and the ``1'' are adjacent in $x$.
\end{lemma}
\begin{proof}
Suppose $x$ contains either $\up{2}31$ or $31\down{2}$.
Specifically let $1\le a<b<c\le n$ be entries in $x$ such that either $bca$ forms the pattern $\up{2}31$ or $cab$ forms the pattern 
$31\down{2}$.
Furthermore choose $(a,b,c)$ among such triples so as to minimize the number of entries occurring between $c$ and $a$.
We claim that $c$ and $a$ are adjacent entries in~$x$.
For the sake of contradiction, suppose there is an entry $d$ between $c$ and $a$.
If $d<b$ then the triple $(d,b,c)$ forms the same pattern as $(a,b,c)$, contradicting our choice of $(a,b,c)$ to minimize the number of entries 
between $c$ and $a$.
On the other hand, if $d>b$, the triple $(a,b,d)$ gives the same contradiction.
The contradiction shows that $c$ and $a$ are adjacent in $x$.
\end{proof}

The inversion set $I(x)$ of a permutation $x\in S_n$ is called {\em $f$-compressed} if the following two conditions hold:
\begin{enumerate}
\item[(i) ] If $i<\upwide{j}<k$ and $(i,k)\in I(x)$ then $(j,k)\in I(x)$, and 
\item[(ii) ] If $i<\downwide{j}<k$ and $(i,k)\in I(x)$ then $(i,j)\in I(x)$.
\end{enumerate}
It is immediate that $I(x)$ is $f$-compressed if and only if $x$ avoids $\up{2}31$ and $31\down{2}$.
This definition is a variant of a definition used in~\cite{Nonpure II} to prove facts about the Tamari lattice.

\begin{lemma}
\label{compressed}
For every $x\in S_n$, there is an element $x'$ such that $I(x')$ is the unique maximal set under containment among all 
$f$-compressed inversion sets $I(w)\subseteq I(x)$.
\end{lemma}
\begin{proof}
By induction on the weak order.
If $I(x)$ is $f$-compressed (as for example when $x$ is the identity permutation) then $x' = x$.
If $I(x)$ is not $f$-compressed, then $x$ contains an instance of either $\up{2}31$ or $31\down{2}$.
By Lemma~\ref{adjacent}, we can take the ``3'' and the ``1'' to be adjacent, so that there is either a $\up{2}31\to\up{2}13$-move or a 
$31\down{2}\to 13\down{2}$-move to some element $z$.
We claim that we can set $x'=z'$.
To prove the claim, let $(i,k)$ be the unique inversion in $I(x)-I(z)$, or in other words let $k$ be the ``3'' and~$i$ be the ``1.''
Let $w\in S_n$ be such that $I(w)$ is $f$-compressed and $I(w)\subseteq I(x)$.
If $z$ is obtained from $x$ by a $\up{2}31\to\up{2}13$-move, then the ``$2$'' is an element $\upwide{j}$ such that $(j,k)\not\in I(x)$.
Since $I(w)\subseteq I(x)$, we have $(j,k)\not\in I(w)$ and since $I(w)$ is $f$-compressed, $(i,k)\not\in I(w)$.
If $z$ is obtained from $x$ by a $31\down{2}\to 13\down{2}$-move, we argue similarly that $(i,k)\not\in I(w)$.
\begin{comment}
If $z$ is obtained from $x$ by a $31\down{2}\to 13\down{2}$-move$, then the ``$\down{2}$'' is an element $\down{j}$ such that 
$(j,i)\not\in I(x)$.
Since $I(w)\subseteq I(x)$, we have $(j,i)\not\in I(w)$ and since $I(w)$ is $f$-compressed, $(k,i)\not\in I(w)$.
\end{comment}
Thus in either case $I(w)\subseteq I(z)$, so by induction $I(w)\subseteq I(z')$.
Since this is true for any $w\in S_n$ such that $I(w)$ is $f$-compressed and $I(w)\subseteq I(x)$, the desired element is $x'=z'$.
\end{proof}

\begin{prop}
\label{projections}
For any permutation $x$, the fiber $\eta^{-1}(\eta(x))$ is the weak order interval $[\pidown x,\piup x]$ for some permutations 
$\pidown x$ and $\piup x$.
Furthermore, $x=\pidown x$ if and only if $x$ avoids $\up{2}31$ and $31\down{2}$, and $x=\piup x$ if and only if $x$ avoids 
$\up{2}13$ and $13\down{2}$.
\end{prop}
\begin{proof}
By Propositions~\ref{cover moves} and~\ref{convex} and Lemma~\ref{adjacent}, an element is minimal in $\eta^{-1}(\eta(x))$ if and 
only if it avoids $\up{2}31$ and $31\down{2}$, or equivalently, if and only if its inversion set if $f$-compressed.
Suppose there is more than one minimal element of $\eta^{-1}(\eta(x))$.
Lemma~\ref{connected} states that $\eta^{-1}(\eta(x))$ is connected, so we can find two distinct minimal elements which are below the same element $y\in\eta^{-1}(\eta(x))$.
Lemma~\ref{compressed} constructs a minimal element $y'$ below $y$ which is (weakly) above any other minimal element, leading to the absurd 
conclusion that $y'$ is minimal and is weakly above two distinct minimal elements.
Thus there is a unique minimal element $\pidown x$.
By symmetry, there is a unique maximal element $\piup x$, and by Proposition~\ref{convex}, $\eta^{-1}(\eta(x))$ is the interval 
$[\pidown x,\piup x]$.
By the symmetry $x\mapsto w_0x$ we have $x=\piup x$ if and only if $x$ avoids $\up{2}13$ and $13\down{2}$.
\end{proof}

We now complete the proof that the fibers of $\eta$ are the congruence classes of a lattice congruence.

\begin{proof}[Proof of Theorem~\ref{A congruence}]
In light of Proposition~\ref{projections}, it remains only to show that $\piup$ and $\pidown$ are order-preserving maps.
Suppose $x\le y$ in $S_n$.
Then $I(\pidown x)$ and $I(\pidown y)$ are both $f$-compressed and $I(\pidown x)\subseteq I(x)\subseteq I(y)$.
By Lemma~\ref{compressed} and Proposition~\ref{projections}, $I(\pidown y)$ is the unique maximal $f$-compressed inversion set
contained in $I(y)$, so $I(\pidown x)\subseteq I(\pidown y)$, or in other words, $\pidown x\le\pidown y$.
The proof that $\piup$ is order-preserving follows by the anti-symmetry $x\mapsto w_0x$.
\end{proof}

\section{Cambrian Lattices of Type A}
\label{Cam A}
By identifying the set of triangulations of $Q$ with the set of permutations $x$ such that $x=\pidown x$, we induce a partial order on the 
triangulations.
The content of Theorem~\ref{A congruence} is that $\eta_f$ is a lattice homomorphism from the weak order onto this partial order, which is 
the quotient lattice $S_n/\Theta_f$.
In this section we show that $S_n/\Theta_f$ is a Cambrian lattice.
We describe the combinatorial realizations of the \mbox{type-A} Cambrian lattices in terms of permutations and of triangulations.
Then we prove the \mbox{type-A} case of Theorem~\ref{ab weak fan}, which states that the Cambrian fans for $S_n$ are combinatorially isomorphic
to the normal fan of the associahedron.
We prove that the \mbox{type-A} Cambrian lattices are sublattices of the weak order on $S_n$, describe the descent map on the \mbox{type-A} Cambrian 
lattices and determine the congruence lattice of a \mbox{type-A} Cambrian lattice.

Recall that a congruence is determined by the set of join-irreducibles it contracts.
A join-irreducible $\gamma$ in $S_n$ is contracted by $\Theta_f$ if and only if it contains the pattern $\up{2}31$ or the pattern $31\down{2}$.
This can be phrased nicely in terms of the subset $A$ associated to $\gamma$.
\begin{lemma}
\label{A ji avoid}
Let $A$ be a subset of $[n]$ corresponding to a join-irreducible element $\gamma$ of the weak order on $S_n$, with $m$ and $M$ as usual.
Then 
\begin{enumerate}
\item[(i) ] $\gamma$ contains $\up{2}31$ if and only if there exists $\up{b}\in A^c\cap(m,M)$.
\item[(ii) ] $\gamma$ contains $31\down{2}$ if and only if there exists $\down{b}\in A\cap(m,M)$.
\end{enumerate}
\end{lemma}
\begin{proof}
Let $\up{b}ca$ be entries in $\gamma$ forming a pattern $\up{2}31$.
Then $b,c\in A^c$ and $a\in A$, and thus $m\le a<b<c\le M$, so that $\up{b}\in A^c\cap(m,M)$.
Conversely, if $\up{b}\in A^c\cap(m,M)$, then $bMm$ forms a $\up{2}31$-pattern in $\gamma$.
The proof for $31\down{2}$ patterns is similar.
\begin{comment}
Let $ca\down{b}$ be entries in $\gamma$ forming a pattern $31\down{2}$.
Then $c\in A^c$ and $a,b\in A$, and thus $m\le a<b<c\le M$, so that $\down{b}\in A\cap(m,M)$.
Conversely, if $\down{b}\in A\cap(m,M)$, then $Mmb$ forms a $31\down{2}$-pattern in $\gamma$.
\end{comment}
\end{proof}

\begin{figure}[ht]
\centerline{\scalebox{.8}{\epsfbox{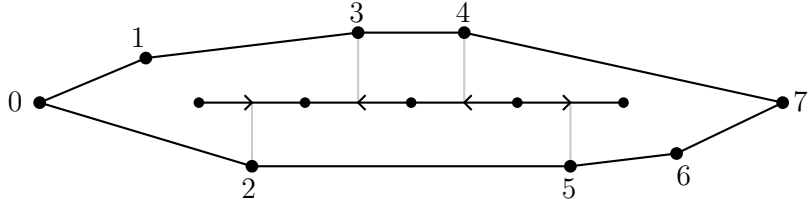}}}   
\caption{An orientation $\G$ and a polygon $Q(\G)$.}
\label{Qorient}
\end{figure}

Orientations $\G$ of the Coxeter diagram for $S_n$ correspond to choices of the linear functional $f$ as follows.
Given $f$, define the orientation $\G_f$ to be $s_b\to s_{b-1}$ for every $\up{b}\in[2,n-1]$ and $s_{b-1}\to s_b$ for every $\down{b}\in [2,n-1]$.
By the reverse process, an orientation $\G$ specifies which indices in $[2,n-1]$ are up or down, and the indices $1$ and $n$ can be arbitrarily 
chosen as up or down indices.
Denote any polygon corresponding to such a choice of up and down vertices by $Q(\G)$.
Figure~\ref{Qorient} illustrates the correspondence between an orientation $\G$ and a choice of up and down vertices.
The polygon shown is a linear scaling of the polygon appearing in Figure~\ref{eta}.

\begin{theorem}
\label{cambrian A}
The congruence $\Theta_f$ is the Cambrian congruence.
In particular the quotient lattice $S_n/\Theta_f$ is the Cambrian lattice $\C(\G_f)$.
\end{theorem}
\begin{proof}
If $\up{b}\in[2,n-1]$, then the join-irreducible permutation 
\[s_{b-1}s_b=12\cdots(b-2)b(b+1)(b-1)(b+2)(b+3)\cdots n\]
 is contracted by $\Theta_f$.
The permutation 
\[s_bs_{b-1}=12\cdots(b-2)(b+1)(b-1)b(b+2)(b+3)\cdots n\]
 is not contracted.
If $\down{b}\in[2,n-1]$, then $s_bs_{b-1}$ is contracted and $s_{b-1}s_b$ is not.
Thus we have verified that $\Theta_f$ contracts the correct set of degree-2 join-irreducibles.
We finish by showing that if $\gamma_2$ is a join-irreducible with $\deg(\gamma_2)>2$, contracted by $\Theta_f$, with associated subset 
$A_2$, then there is a join irreducible $\gamma_1$, also contracted by $\Theta_f$, whose associated subset $A_1$ has $A_1\to A_2$.
This will complete the proof since $\Irr(\Con(S_n))$ is dually ranked by degree.

We use Lemma~\ref{A ji avoid}.
Suppose either $\up{b}\in A_2^c\cap(m_2,M_2)$ or $\down{b}\in A_2\cap(m_2,M_2)$.
Since $\deg(\gamma_2)>2$, or in other words $M_2-m_2>2$, we have either $b+1\in(m_2,M_2)$ or $b-1\in(m_2,M_2)$.
If $b+1\in(m_2,M_2)$, let $A_1=\left(A_2\cap[1,b]\right)\cup[b+2,n]$.  
Thus $m_1=m_2$, $M_1=b+1<M_2$ and $A_1\cap[1,M_1)=A_2\cap[1,M_1)$.
Also, $b\in A_1$ if and only $b\in A_2$, so either $\up{b}\in A_1^c\cap(m_1,M_1)$ or $\down{b}\in A_1\cap(m_1,M_1)$.
If $b-1\in(m_2,M_2)$, let $A_1=\set{b-1}\cup\left(A_2\cap[b,n]\right)$.  
Thus $m_1=b-1>m_2$, $M_1=M_2$ and $A_1\cap(m_1,n]=A_2\cap(m_1,n]$.
Again, $b\in A_1$ if and only $b\in A_2$, so either $\up{b}\in A_1^c\cap(m_1,M_1)$ or $\down{b}\in A_1\cap(m_1,M_1)$.
\end{proof}

\begin{figure}[ht]
\centerline{\scalebox{.6}{\epsfbox{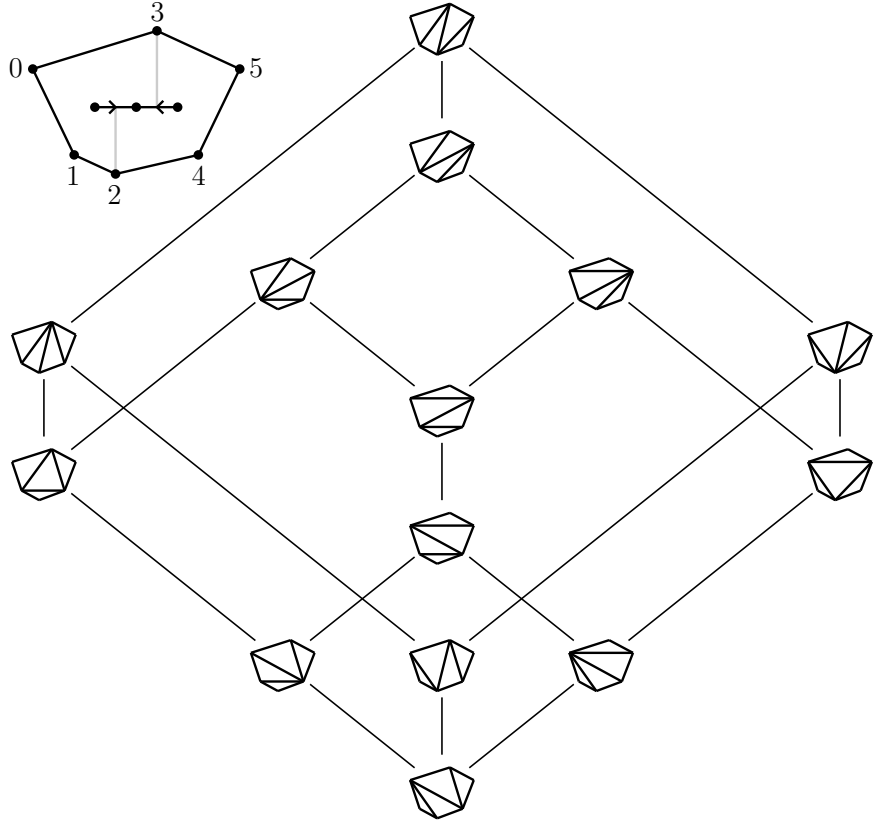}}}   
\caption{A Cambrian lattice.}
\label{A3camb_tri}
\end{figure}

As an immediate consequence of Theorem~\ref{cambrian A} and the results of the preceding sections, we have the following combinatorial
characterizations of the Cambrian lattices of type~A.
Theorem~\ref{tri A} is illustrated in Figure~\ref{A3camb_tri}, a continuation of the $S_4$ example of Figure~\ref{A3tri}.
\begin{theorem}
\label{tri A}
For $\G$ an orientation of the Coxeter diagram for $S_n$, the Cambrian lattice $\C(\G)$ is isomorphic to the partial order on triangulations 
of $Q(\G)$  whose cover relations are diagonal flips, where going up in the cover relation corresponds to increasing the slope of 
the diagonal.
\end{theorem}
\begin{theorem}
\label{per A}
For $\G$ an orientation of the Coxeter diagram for $S_n$, the Cambrian lattice $\C(\G)$ is isomorphic to the subposet of the weak order on 
$S_n$ consisting of permutations avoiding both $\up{2}31$ and $31\down{2}$, where up and down indices are determined by the vertices of 
$Q(\G)$.
\end{theorem}

The \mbox{type-A} case of Theorem~\ref{ab weak fan} is also a consequence of Theorem~\ref{cambrian A}, as we now explain.
The iterated fiber polytope $\Sigma(\Delta\to Q\to I)$ is combinatorially isomorphic to the permutohedron, and the fiber polytope 
$\Sigma(\Delta\to Q)$ is combinatorially isomorphic to the associahedron.
Thus the normal fan of $\Sigma(\Delta\to Q\to I)$ is combinatorially isomorphic to the normal fan $\F$ of the permutohedron.
The maximal cones of $\F_{(\Theta_f)}$ are the unions over $\Theta_f$-classes of the maximal cones of $\F$.
But each $\Theta_f$-class corresponds to the set of maximal normal cones to $\Sigma(\Delta\to Q\to I)$ contained in some maximal normal 
cone to $\Sigma(\Delta\to Q)$.
Thus $\F_{(\Theta_f)}$ is also combinatorially isomorphic to the associahedron.

The following theorem strengthens Theorem~\ref{per A}.
The proof is similar to the proof given in~\cite{Nonpure II} for the Tamari lattice.
\begin{theorem}
\label{sub A}
For any orientation $\G$ of the Coxeter diagram associated to $S_n$, the Cambrian lattice $\C(\G)$ is a sublattice of the weak order on 
$S_n$.
\end{theorem}
\begin{proof}
Choose an $f$ such that $\G=\G_f$, and let $\piup$ and $\pidown$ correspond to the congruence $\Theta_f$.
Let $x,y\in S_n$ have $\pidown x=x$ and $\pidown y=y$.
Since $I(x)$ is compressed and $I(x)\subseteq I(x\join y)$, by Lemma~\ref{compressed} we have 
$I(x)\subseteq I(\pidown(x\join y))$, or in other words $x\le \pidown(x\join y)$.
Similarly, $y\le \pidown(x\join y)$, so $x\join y\le\pidown (x\join y)$ and therefore $x\join y=\pidown (x\join y)$.

Let $z:=x\meet y$.
We need to show that $\pidown z=z$.
Suppose for the sake of contradiction that $\pidown z<z$, or in other words that $z$ contains either $\up{2}31$ or $31\down{2}$.
If $z$ contains the pattern $\up{2}31$, let $1\le a<b<c\le n$ be entries in $z$ such that $bca$ forms the pattern $\up{2}31$.
Furthermore choose $(a,b,c)$ among such triples so as to minimize the number of entries occurring between $b$ and $c$ in $z$.
Since $(a,c)\in I(z)$, we have $(a,c)\in I(x)$, and since $I(x)$ is $f$-compressed, $(b,c)\in I(x)$.
Similarly $(b,c)\in I(y)$.

Let $d$ be the entry of $z$ immediately preceding $c$.
Then by our choice of $(a,b,c)$, either $d=b$ or $d<b$, and in either case, $d<c$.
Let $w$ be the permutation obtained from $z$ by transposing $d$ and $c$, so that in particular $z\covered w$.
Also, $I(w)=I(z)\cup\set{(d,c)}$.
If $d<b$, then $(d,b)\in I(z)\subseteq I(x)$, and because $(b,c)\in I(x)$, by transitivity we have $(d,c)\in I(x)$.
If $d=b$, then $(d,c)=(b,c)\in I(x)$, so in either case we have $I(w)\subseteq I(x)$.
The same reasoning shows $I(w)\subseteq I(y)$, and this contradicts the fact that $z=x\meet y$, thus proving that $\pidown z=z$.

If $z$ contains the pattern $31\down{2}$ we argue symmetrically by the symmetry $x\mapsto w_0xw_0$ to obtain the same contradiction.
\begin{comment}
We find an instance $cab$ of this pattern so as to minimize the number of entries occurring between $a$ and $b$.
By considering the entry $d$ immediately following $a$ and reasoning similarly to the previous case, we arrive at the same contradiction.
Since $(c,a)\in I(z)$, we have $(c,a)\in I(x)$, and since $I(x)$ is $f$-compressed, $(b,a)\in I(x)$.
Similarly $(b,a)\in I(y)$.
Let $d$ be the entry of $z$ immediately following $a$.
Then either $d=b$ or $d>b$, and in any case $d>a$.
Let $w$ be the permutation obtained from $z$ by transposing $a$ and $d$, so that $z\covered w$ and $I(w)=I(z)\cup\set{(d,a)}$.
If $d>b$, then $(d,b)\in I(z)\subseteq I(x)$, and because $(b,a)\in I(x)$, by transitivity we have $(d,a)\in I(x)$.
If $d=b$, then $(d,a)=(b,a)\in I(x)$, so in either case we have $I(w)\subseteq I(x)$.
The same reasoning shows $I(w)\subseteq I(y)$, and this contradicts the fact that $z=x\meet y$.
\end{comment}
\end{proof}

Recall from Section~\ref{cong} that if $\eta$ is the canonical homomorphism from $W$ to a Cambrian lattice $\C(\G)$ then the descent 
map $\delta$ on permutations factors through~$\eta$.
This can be seen in the special case $W=S_n$ by noting that the moves leading from $x$ to $\pidown x$ never destroy left descents.
Thus we can define the descent set of a triangulation $T$ to be the left descent set of any $x\in S_n$ with $\eta(x)=T$.
The following proposition follows immediately from the characterization of $\eta$ in terms of paths $\lambda$.
As one might expect, determining whether $(a,a+1)$ is a descent depends on whether $a$ and $a+1$ are up indices or down indices.

\begin{prop}
\label{A descents}
Let $T$ be a triangulation of $Q$.
The descent set of $T$ is described as follows.  
For $a\in[n-1]$,
\begin{enumerate}
\item[(i) ]$(\down{a},\down{a+1})$ is a descent if and only if there is a diagonal $a$---$i$ in $T$ with $i>a+1$, which is if and only if
there is no diagonal $i$---$(a+1)$ with $i<a$.
\item[(ii) ]$(\down{a},\up{a+1})$ is a descent if and only if $a$---$(a+1)$ is a diagonal in $T$.
\item[(iii) ]$(\up{a},\up{a+1})$ is a descent if and only if there is no diagonal $a$---$i$ in $T$ 
with $i>a+1$, which is if and only if there is a diagonal $i$---$(a+1)$ with $i<a$.
\item[(iv) ]$(\up{a},\down{a+1})$ is a descent if and only if $a$---$(a+1)$ is not a diagonal in $T$.
\end{enumerate}
\end{prop}
We conclude the section by determining the congruence lattice of a Cambrian lattice of type~A\@.
\begin{prop}
\label{A one each}
For any orientation $\G$ of the Coxeter diagram for $S_n$, given a reflection $(m,M)$, there is a unique join-irreducible $\gamma$ which is not 
contracted by $\Theta(\G)$ and whose associated left reflection is $(m,M)$.
The subset associated to $\gamma$ is $A=\set{m}\cup\up{(m,M)}\cup(M,n]$.
\end{prop}
\begin{proof}
Let $\gamma$ be a join-irreducible whose associated left reflection is $(m,M)$ and let $A$ be the subset corresponding to a $\gamma$.
By Lemma~\ref{A ji avoid}, if $\gamma$ avoids $\up{2}31$ then $A^c\cap\up{(m,M)}=\emptyset$, and if $\gamma$ avoids $31\down{2}$ 
then $A\cap\down{(m,M)}=\emptyset$.
Thus if $\gamma$ avoids both then $A=\set{m}\cup\up{(m,M)}\cup(M,n]$.
\end{proof}
\begin{prop}
\label{A camb cong}
Suppose $\gamma_1$ and $\gamma_2$ are join-irreducibles in $S_n$ not contracted by $\Theta(\G)$, with associated subsets $A_1$ and 
$A_2$ and associated left reflections $(m_1,M_1)$ and $(m_2,M_2)$.
Then $\gamma_1\to \gamma_2$ in the directed graph of Theorem~\ref{A shard} if and only if either $m_2=m_1<M_1<M_2$ or 
$m_1<m_2<M_2=M_1$.
\end{prop}
\begin{proof}
We apply Theorem~\ref{A shard}.
The ``only if'' statement is immediate, and if either condition above holds then by Proposition~\ref{A one each} we have
$A_1\cap(m_1,M_1)=\up{(m_1,M_1)}=A_2\cap(m_1,M_1)$.
\end{proof}

Readers familiar with root systems will notice that Proposition~\ref{A camb cong} implies that for any orientation $\G$ of the diagram 
for $S_n$, the poset $\Irr(\Con(\C(\G)))$ is isomorphic to the root poset of $S_n$. 
(This is not true for general root systems.)
In particular, the number of congruences of a Cambrian lattice of type~A is the Catalan number. 
This extends to all Cambrian lattices of type~A the result of Geyer~\cite{Geyer} on the number of congruences of the Tamari lattice.
Geyer also showed that the number of cover relations in the congruence lattice is equal to the number of edges in the associahedron, and that 
result also extends to all Cambrian lattices of type~A\@.

\section{Cambrian Lattices of Type B}
\label{Cam B}
In this section we review an equivariant version, due to Reiner~\cite{Equivariant}, of the fiber polytope construction described previously, 
leading to a map from the $B_n$-permutohedron to the $B_n$-associahedron.
We use this map to construct the Cambrian lattices of type~B\@.
The equivariant fiber-polytope construction produces combinatorial realizations which are related to the \mbox{type-A} case by the standard ``folding'' 
construction.
However, the fact that these are combinatorial realizations of the Cambrian lattices does not follow in any obvious way from the \mbox{type-A} proof 
by folding, but must be argued separately.
We describe the combinatorics of the \mbox{type-B} Cambrian lattices, show that they are sublattices of the weak order on $B_n$ and describe the
descent map on these lattices.
Finally, we identify the \mbox{type-B} Tamari lattices and describe them in terms of signed pattern avoidance.

Returning to the general case of two polytopes $P$ and $Q$ with a linear surjection $P\mapname{\varphi}Q$, we now suppose that some group $G$
acts on both $P$ and $Q$ by symmetries such that $\varphi$ is a $G$-equivariant map.
Then $G$ also acts on the fiber polytope \mbox{$\Sigma(P\mapname{\varphi}Q)$} by symmetries, and the {\em equivariant fiber polytope} 
$\Sigma^G(P\mapname{\varphi}Q)$ can be defined as the set of points in the fiber polytope which are fixed by this $G$-action.
Given a group $G$ acting by symmetries on polytopes $P$, $Q$ and $R$ and a tower $P\mapname{\varphi}Q\mapname{\rho}R$ of surjective 
$G$-equivariant linear maps, the map $\varphi$ from $\Sigma(Q\mapname{\rho}R)$ to $\Sigma(P\mapname{\rho\circ\varphi}R)$ is 
$G$-equivariant.
Thus there is an {\em iterated equivariant fiber polytope} $\Sigma^G(P\mapname{\varphi}Q\mapname{\rho}R)$ defined as
$\Sigma^G\left(\Sigma(P\mapname{\rho\circ\varphi}R)\mapname{\varphi}\Sigma(Q\mapname{\rho}R)\right)$.
By restriction of the non-equivariant result to the subspace of fixed points, we have that the normal fan of 
$\Sigma^G(P\mapname{\varphi}Q\mapname{\rho}R)$ refines that of $\Sigma^G(P\mapname{\varphi}Q)$.

We consider the two-element group $C_2$ acting on the polytopes of Section~\ref{A assoc}.
Define a tower of surjective linear maps of polytopes 
\[\Delta^{2n+1}\mapname{\varphi} Q_{2n+2}\mapname{\rho} I,\]
where $I$ is again the interval, $Q_{2n+2}$ is a centrally symmetric $(2n+2)$-gon and $\Delta^{2n+1}$ is a simplex with vertices 
$ e_{-n-1}, e_{-n},\ldots, e_{-1}, e_1, e_2,\ldots, e_{n+1}$.
Once again $v_i=\varphi( e_i)$ and $a_i=\rho(v_i)$ with $a_i<a_j$ whenever $i<j$.
Require that the central symmetry of $Q_{2n+2}$ map $v_i\mapsto v_{-i}$ and that the central symmetry of $I$ map $a_i\mapsto a_{-i}$.
Let the non-trivial element of $C_2$ act respectively on $Q$ and $I$ by these symmetries and on $\Delta$ by $ e_i\mapsto  e_{-i}$.
Let $f$ be a non-trivial linear functional on $\ker\rho$, and use the notation $\upwide{i}$ and $\downwide{i}$ as before.
Then we have $f_i=-f_{-i}$ or in other words,~$i$ is $\upwide{i}$ if and only if $-i$ is $\downwide{-i}$.

The iterated fiber polytope $\Sigma(\Delta\mapname{\varphi}Q\mapname{\rho}I)$ is the $S_{2n}$-permutohedron.
The $C_2$-equivariant version is the $B_n$-permutohedron, with vertices given by signed permutations.
This is because $C_2$ acts on permutations of $\pm[n]$ by the automorphism $w\mapsto w_0ww_0$, which reverses the order of the entries in 
the permutation and reverses the sign of each entry.
Thus the fixed permutations are exactly the signed permutations. 
The fiber polytope $\Sigma(\Delta\mapname{\varphi} Q)$ is the $S_{2n}$-associahedron, and the $C_2$-equivariant version is
the $B_n$-associahedron, whose vertices are triangulations of $Q_{n+2}$ which are fixed by the $C_2$-action generated by 
$v_i\mapsto v_{-i}$.
For a fixed choice of $f$, there is a map $\eta_f^B$ from the vertices of $\Sigma^{C_2}(\Delta\to Q\to I)$ to the vertices of 
$\Sigma^{C_2}(\Delta\to Q)$.
This is the restriction to signed permutations of the map $\eta_f$ from permutations of $\pm[n]$ to triangulations of $Q_{2n+2}$.
The map $\eta_f^B$ on signed permutations can be characterized in terms of paths $\lambda$ in exactly the same way as for unsigned 
permutations.
By Propositions~\ref{auto sublattice} and~\ref{sub congruence}, we have the following theorem.

\begin{theorem}
\label{B congruence}
The fibers of $\eta_f^B$ are the congruence classes of a lattice congruence $\Theta_f^B$ of the weak order on $B_n$.
\end{theorem}
In particular, $\eta_f^B$ induces a partial order on the centrally symmetric triangulations of $Q$, which are the vertices of the \mbox{type-B} 
associahedron.
This partial order is a lattice homomorphic image of the weak order.
Identify these vertices with the signed permutations $x$ such that  $x=\pidown x$.

We construct an orientation $\G^B_f$ of the Coxeter diagram for $B_n$ by directing \mbox{$s_b\to s_{b-1}$} for every $\up{b}\in[n-1]$ and 
$s_{b-1}\to s_b$ for every $\down{b}\in [n-1]$.
Given an orientation~$\G$, one can conversely construct a centrally symmetric polygon $Q(\G)$.
A join-irreducible $\gamma$ of $B_n$, with corresponding signed subset $A$ is contracted by $\Theta_f^B$ if and only if the full notation for 
$\gamma$, thought of as a permutation of $\pm[n]$, contains the pattern $\up{2}31$ or equivalently by the symmetry of signed permutations, 
the pattern $31\down{2}$.
Here we don't want signed pattern avoidance, but rather the variant of pattern avoidance defined at the beginning of Section~\ref{eta A}.
Since $\gamma$ consists of the elements of $-A$ in ascending order, followed by $(\pm A)^c$ in ascending order, then the elements of $A$ in 
ascending order, we can phrase this nicely in terms of $A$.

\begin{lemma}
\label{B ji avoid}
Let $A$ be a signed subset of $[n]$ corresponding to a join-irreducible element $\gamma$ of the weak order on $B_n$, with $m$ and $M$ as above.
Then $\gamma$ contains $\up{2}31$ or equivalently $31\down{2}$ if and only if one of more of the following holds:
\begin{enumerate}
\item[(i) ] There exists $\up{b}\in A^c\cap(m,M)$, and/or
\item[(ii) ] There exists $\down{b}\in A\cap(m,M)$.
\end{enumerate}
\end{lemma}
\begin{proof}
If $\up{b}ca$ are entries in $\gamma$ forming a pattern $\up{2}31$, then one of the following holds:
\begin{itemize}
\item[(A) ]$b,c\in(\pm A)^c$ and $a\in A$,
\item[(B) ]$b\in -A$, $c\in (\pm A)^c$ and $a\in A$,
\item[(C) ]$b,c\in -A$ and $a\in A$, or
\item[(D) ]$b,c\in -A$ and $a\in (\pm A)^c$.
\end{itemize}
In cases (A) and (B) we have $m\le a<b<c\le M$, so (i) holds.
In case (C) we have $m\le a<b<c\le -m$ so that $\up{b}\in(m,-m)$, and therefore $\down{-b}\in(m,-m)$.
If $\down{-b}\in(m,M)$ this is (ii).
Otherwise, we have $\down{-b}\in[M,-m)$, so that $\up{b}\in(m,-M]\subseteq(m,M)$ and $b\in(-A)\subseteq A^c$, so that (i) holds.
In case (D) we have $-M\le a<b<c\le -m$, so $\down{-b}\in(m,M)$, and this is (ii).

Conversely, if (i) holds then $\up{b}Mm$ is the pattern $\up{2}31$ in $\gamma$ and if (ii) holds then $Mm\down{b}$ is the pattern
$31\down{2}$ in $\gamma$.
\end{proof}

\begin{theorem}
\label{cambrian B}
The quotient lattice $B_n/\Theta^B_f$ is the Cambrian lattice $\C(\G^B_f)$.
\end{theorem}
\begin{proof}
The following facts are easily checked using Lemma~\ref{B ji avoid}:
The join-irreducibles $s_0s_1$ and $s_0s_1s_0$ are contracted by $\Theta_f^B$ if and only if $1$ is $\up{1}$.
The join-irreducibles $s_1s_0$ and $s_1s_0s_1$ are contracted if and only if $1$ is $\down{1}$.
For $b\in[2,n-1]$, the join-irreducible $s_bs_{b-1}$ is contracted if and only if $b$ is $\down{b}$, and
the join-irreducible $s_{b-1}s_b$ is contracted if and only if $b$ is $\up{b}$.
Thus the set of contracted degree-2 join irreducibles is correct, and it remains to show that $\Theta^B_f$ is generated in degree 2.

Suppose that $\gamma_2$ is a join-irreducible with $\deg(\gamma_2)>2$, contracted by $\Theta^B_f$, with associated signed subset $A_2$.
Thus either condition (i) or condition (ii) of Lemma~\ref{B ji avoid} holds.
For the remainder of the proof, $b$ stands for a value such that condition (i) or (ii) holds for $A_2$.
We will show that there is a join irreducible $\gamma_1$, also contracted by $\Theta^B_f$, whose associated signed subset $A_1$ has $A_1\to A_2$.
In many cases we will construct an $A_1$ with $A_1\to A_2$ such that (r1) holds, and such that $b$ is in $(m_1,M_1)$.
That finishes the proof for that case because condition (r1) assures that $b\in A_1$ if and only if $b\in A_2$, and thus condition (i) or (ii) holds
for $A_1$ as well.

First, notice that the case $m_2>0$ can be dealt with exactly as in the proof of Theorem~\ref{cambrian A}, so from now on we assume $m_2<0$.
Next consider the case where $m_2=-M_2$.
Since in this case $A_2^c=-A_2$, we have $\up{b}\in A_2^c\cap(m_2,M_2)$ if and only if $\down{-b}\in A_2\cap(m_2,M_2)$.
Thus we may as well assume $b>0$.
By Lemma~\ref{B degrees}, we must have $M_2>2$.
If $b\neq 1$ then let $A_1=\set{1}\cup(A_2\cap(1,M_2))\cup(M_2,n]$.
Then $-m_2=M_2=M_1>m_1=1>0$, so (q2) holds, and $A_1$ was defined so that (r1) holds as well.
Thus $A_1\to A_2$, and since $b\in(m_1,M_1)$ we are finished.
If $b=1$ then let $A_1=(A-\set{m,-m-1})\cup\set{m+1,-m}$.
Then $-m_1=M_1=-m_2-1<M_2=-m_2$, so (q1) holds, and (r1) holds by construction.
Since $M_2>2$, we have $b\in(m_1,M_1)$ in this case as well.

From now on assume $M_2\neq-m_2$ and $m_2<0$.
If 
\[(\pm A_2)^c\cap(m_2,M_2)\cap(-M_2,-m_2)\neq\emptyset,\]
then since this set is sign-symmetric, there is a $\up{b}\in (\pm A_2)^c\cap(m_2,M_2)\cap(-M_2,-m_2)$.
Thus condition (i) holds with this $b$.
Notice that (f2\,:\,$-b$) holds because $-b\in A_2^c-\set{-M_2,-m_2}$ and $b\in A_2^c$.
If $b>0$ then let $A_1=(A_2-[-n,-b))\cup\set{-b}\cup(M_2,n]$.
Then $M_2=M_1>m_1=-b>m_2\neq-M_2$, so (q3) holds and by construction (r1) holds.
Since we have (r1) and $b\in(m_1,M_1)$, we are finished.
If $b<0$, let $A_1=A_2\cup(-A_2^c\cap(-b,n])$.
Thus $M_2>M_1=-b>m_1=m_2\neq-M_2$, so condition (q4) holds, (r1) holds by construction, and $b\in(m_1,M_1)$.

Thus from now on, we continue to assume $M_2\neq-m_2$ and $m_2<0$ and begin also to assume 
$(\pm A_2)^c\cap(m_2,M_2)\cap(-M_2,-m_2)=\emptyset$, which we refer to as ``the emptiness assumption.''
The proof can now be finished by considering three cases, depending on where $b$ lies in relation to $m_2$ and $M_2$.
\begin{enumerate}
\item[]
\begin{enumerate}
\case{1}{$b\in(m_2,M_2)\cap(-M_2,-m_2)$}
Let $\mu=\max\set{m_2,-M_2}$.
Since $b\in(m_2,M_2)\cap(-M_2,-m_2)$, we must have $\mu\le -2$.
Let $A_1=(A_2-[-n,\mu))\cup\set{\mu}\cup(-\mu,n]$.
Then $m_1=\mu$ and by the emptiness assumption and by construction we have $|A_1|=n$, so $M_1=-m_1$.
Also by construction, either (q3) or (q4) holds, and (r1) holds.
By the emptiness assumption, either (f3\,:\,$m_1$) or (f3\,:\,$M_1$) holds as required.
\case{2}{$m_2<-m_2\le b<M_2$}
If $m_2<-1$ then by the emptiness assumption, $|A_2\cap\set{1,-1}|=1$, so let $\mu$ be the element of that intersection.
Let $A_1=(A_2\cap[\mu,n])$, so that $m_1=\mu$ and $M_1=M_2$.
By construction (q3) and (r1) hold, and (f1\,:\,$m_1$) holds because $m_1\in A_2$.
If on the other hand $m_2=-1$ and $b>1$, then let $A_1=A_2-\set{-1}\cup\set{1}$, so that (q3) and (r1) hold.
Also (f3\,:\,1) holds because $1=-m_2$ and by the emptiness assumption.
Finally, if $m_2=-1$ and $b=1$ then by Lemma~\ref{B degrees}, since $\deg(\gamma_2)>2$ we have $M_2\ge 3$.
Let $A_1=A_2-\set{M_2-1}\cup\set{M_2}$, so that $m_1=m_2$, $M_1=M_2-1$ and therefore (q4) and (r1) hold.
We have (f\,:\,$M_1$) because either $M_1\in A_2$ or $M_1\in A_2^c\cap(-m_2,n]$ and $-M_1\not\in(m_2,M_2)$.
\case{3}{$m_2<b\le-M_2<M_2$}
If $m_2<b-1$ then let $A_1=(A_2\cup\set{m_2+1,-m_2})-\set{m_2,-m_2-1}$.
Then $m_1=m_2+1$ and $M_1=M_2$ so (q3) holds, and by construction (r1) holds.
If $m_1\in A_2$ then (f1\,:\,$m_1$) holds.
If $m_1\in A^c_2$ then since $m_1<-M_2$ we have $-m_1\not\in(m_2,M_2)$, so (f2\,:\,$m_1$) holds.

If $m_2=b-1$ then by Lemma~\ref{B degrees} since $\deg(\gamma_2)>2$, we have $b\le-2$.
Let $\mu$ be the maximal element of $A_2^c\cap\set{-1,1}$ and let $M=-m_2$.
Now let $A_1=\set{\mu}\cup(-A_2^c\cap(\mu,M))\cup(M,n]$.
This is a signed subset because by the emptiness assumption $-A_2^c\cap(\mu,M)=A_2\cap(\mu,M)$.
Then $m_1=\mu$ and $M_1=M$.
If $M_2>1$ then $m_1=\mu>-M_2$ and if $M_2=1$ then $\mu=1>-M_2$.
Thus (q5) holds and (r2) holds by construction.
If $M_2>1$ then by the emptiness assumption $-\mu\in A_2$, so (f1\,:\,$-m_1$) holds.
If $M_2=1$ then $\mu=1$ so $-\mu=-M_2$ and (f3\,:\,$-m_1$) holds by the emptiness assumption.
Now since $b=m_2+1\le-2$ and $m_1=\pm 1$ we have $-b\in(1,-m_2)\subseteq(m_1,M_1)$.
If $\up{b}\in A_2^c$ then $\down{-b}\in -A_2^c$, and therefore $\down{-b}\in A_1$.
If $\down{b}\in A_2$ then $\up{-b}\in-A_2$, and therefore $\up{-b}\in A_1^c$.
Thus $A_1\to A_2$ and $A_1$ satisfies condition (i) and/or condition (ii) of Lemma~\ref{B ji avoid}.
\end{enumerate}
\end{enumerate}
\end{proof}

The \mbox{type-B} case of Theorem~\ref{ab weak fan} can now be argued exactly as in type~A, changing all references to \mbox{type-A} objects to their
\mbox{type-B} analogs.
The following theorems follow from Theorems~\ref{tri A}, \ref{per A} and~\ref{cambrian B} and the fact that $\C(\G)$ is the set of elements 
of the \mbox{type-A} construction which are fixed by the action of $C_2$.
Say a diagonal of $Q$ is a {\em diameter} if its endpoints are related by central symmetry.

\begin{theorem}
\label{tri B}
For $\G$ an orientation of the Coxeter diagram for $B_n$, the Cambrian lattice $\C(\G)$ is isomorphic to the partial order on centrally 
symmetric triangulations of $Q(\G)$ whose cover relations are diameter flips or symmetric pairs of diagonal flips, where going up in the 
cover relation corresponds to increasing the slope of the diameter or of each diagonal.
\end{theorem}
\begin{theorem}
\label{per B}
For $\G$ an orientation of the Coxeter diagram for $B_n$, the Cambrian lattice $\C(\G)$ is isomorphic to the subposet of the weak order on 
$B_n$ consisting of signed permutations whose full notations, thought of as permutations of $\pm[n]$, avoid $\up{2}31$, where up and down 
indices are determined by the vertices of $Q(\G)$.
\end{theorem}

Theorem~\ref{sub A} and Proposition~\ref{auto sublattice} show that $\C(\G)$ is the intersection of two sublattices of the weak order 
on permutations.
Thus we have the following.
\begin{theorem}
\label{sub B}
For any orientation $\G$ of the Coxeter diagram associated to $B_n$, the Cambrian lattice $\C(\G)$ is a sublattice of the weak order on $B_n$.
\end{theorem}

The descent set of a centrally symmetric triangulation of $Q$ is easily described.
\begin{prop}
\label{B descents}
Let $T$ be a centrally symmetric triangulation of $Q$.
For $a\in[n-1]$, one can determine whether $(a,a+1)$ is a descent of $T$ exactly as in Theorem \ref{A descents}.
In addition, 
\begin{enumerate}
\item[(v) ]$(\down{-1},\up{1})$ is a descent if and only if $(-1)$---$1$ is an edge in $T$.
\item[(vi) ]$(\up{-1},\down{1})$ is a descent if and only if $(-1)$---$1$ is not an edge in $T$.\mbox{\huge\protect\phantom{X}}
\end{enumerate}
\end{prop}

\begin{figure}[ht]
\centerline{\scalebox{0.8}{\epsfbox{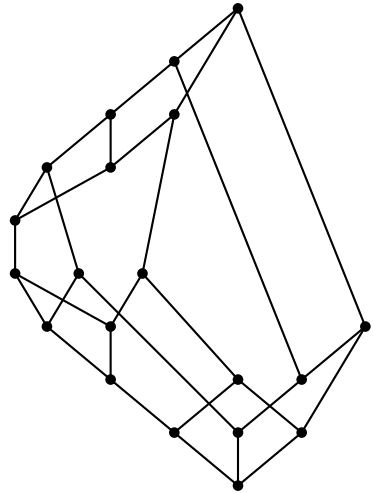}}}   
\caption{A $B_3$ Tamari lattice.}
\label{B3tam}
\end{figure}

In the special cases where $\G$ is directed linearly from one endpoint of $G$ to the other, we call $\C(\G)$ a \mbox{type-B} Tamari lattice.
The justification for the name comes from the fact that, in analogy to type~A, these are the unique Cambrian lattices which can be defined via
signed pattern avoidance (see Section~\ref{weak}), without reference to up indices and down indices.
Recall that by Theorem~\ref{iso} these two orientations define anti-isomorphic lattices.
Figure~\ref{B3tam} shows the $B_3$ Tamari lattice arising from the orientation $s_0\to s_1\to s_2$.
To avoid clutter, we do not label the elements by signed permutations.
The elements covering the minimal element, from left to right, are $s_0$, $s_1$ and $s_2$.

\begin{theorem}
\label{B Tamari}
One of the \mbox{type-B} Tamari lattices is the sublattice of the weak order on signed permutations consisting of signed permutations avoiding the 
signed patterns $-2$$-1$, $2$$-1$, $-231$, $-12$$-3$, $12$$-3$ and $231$.
The other is the sublattice consisting of signed permutations avoiding $-21$, $1$$-2$, $-2$$-1$$-3$, $-13$$-2$, $3$$-12$ and~$312$.
\end{theorem}

\begin{proof}
Any Cambrian lattice of type~B is the sublattice of the weak order on signed permutations consisting of signed permutations whose full notation,
thought of as a permutation of $\pm[n]$, avoids $\up{2}31$.
To relate this condition to signed pattern avoidance in the abbreviated notation, we use the symbol ``$|$'' to represent the position of the 
center of the full notation for a permutation. 
So for example, we say that $x$ contains 4$-6$$|$3, meaning that $x_i=4$, $x_j=-6$ and $x_k=3$ for some  $i<j<0<k$.
We also refer to entries occurring ``on the left'' or ``on the right'' of $x$ to indicate whether they are left or right of the center.
In every case, the phrase ``signed pattern'' is short for ``signed pattern in the abbreviated notation for $x$.'' 

One of the \mbox{type-B} Tamari lattices is characterized by the condition that $\upwide{i}$ if and only if $i>0$.
We will show that in this case, containing $\up{2}31$ is equivalent to containing one or more of the six signed permutations listed.
Suppose the full notation for a signed permutation $x$ contains a subsequence $bca$ such that $a<b<c$ and $0<b$.
If $a=-b$, then $x$ contains either $b(-c)|ca$ or $bc|(-c)a$, so that either $ca$ forms the signed pattern 2$-1$ or $(-c)a$ forms 
the signed pattern $-2$$-1$.
If $a=-c$, then $x$ contains $bc|a(-b)$, so that $a(-b)$ forms the signed pattern -2-1.
Now suppose that $a\neq-b$ and $a\neq-c$.
If $b$ is on the right, there are several cases:
If $0<a<b$ then $bca$ forms the signed pattern 231, if $-c<a<0$ then $ca$ forms the signed pattern 2$-1$ and if $a<-c$ then $bca$ forms the signed pattern 
12$-3$.
If $b$ and $c$ are both on the left, then $(-c)(-b)$ forms the signed pattern $-2$$-1$.
Thus it remains to consider the case where $b$ is on the left and $c$ is on the right.
If $-b$ is to the right of $c$ then $c(-b)$ forms the signed pattern 2$-1$, so we assume that $b|(-b)ca$ occurs in $x$.
If $-c<a<0$ then $ca$ forms the signed pattern 2-1, and otherwise $(-b)ca$ forms either the signed pattern $-231$ or the signed pattern $-12$$-3$.
Conversely, if $x$ contains any of the six signed patterns listed, then the argument is easily reversed to identify an occurrence of $\up{2}31$ in the 
full notation for $x$.

For the other \mbox{type-B} Tamari lattice, we have $\upwide{i}$ if and only if $i<0$.
Suppose that $x$ contains a subsequence $bca$ such that $a<b<c$ and $b<0$.
If $c=-a$ then $x$ contains $bc|a(-b)$, so $a(-b)$ forms the signed pattern $-21$.
If $c=-b$ then $x$ contains $b|ca$ and $ca$ forms the signed pattern $1$$-2$.
Now suppose that $c\neq-a$ and $c\neq-b$.
Suppose $b$ is on the right.
If $b<c<0$ then $bca$ forms the signed pattern $-2$$-1$$-3$, if $0<c<-a$ then $ca$ forms the signed pattern 1$-2$, and if $-a<c$ then $bca$ forms the
signed pattern $-13$$-2$.
If $b$ is on the left and $a$ is on the right, then either $(-b)a$ forms the signed pattern 1$-2$ or $a(-b)$ forms the signed pattern $-2$1.
Finally, suppose $a$ is on the left, so that $|(-a)(-c)(-b)$ occurs in $x$.
If $b<c<0$ then $(-a)(-c)(-b)$ forms the signed pattern 312, if $0<c<-b$ then $(-a)(-c)(-b)$ forms the signed pattern 3$-12$, and if $-b<c$ then 
$(-c)(-b)$ forms the signed pattern $-21$.
Again, the argument is easily reversed.
\end{proof}

\section{Clusters}
\label{cluster}
In this section we give a brief overview of the construction of clusters, and prove Theorem~\ref{cluster refine}.
We conclude by summarizing the combinatorial realization~\cite{ga} of \mbox{type-A} clusters as triangulations.
We assume basic definitions and results on root systems, for which the reader is referred to~\cite{Bourbaki,Humphreys}.
For more details on clusters, see~\cite{ga}.

Let $\Phi$ be a finite crystallographic root system spanning the vector space $\reals^n$ and let $\Pi$ be the simple roots of $\Phi$.
Let $(W,S)$ be the corresponding Coxeter system, with $|S|=n$.
The Coxeter diagram of any finite Coxeter group is a labeled bipartite graph.
Let $I_+$ and $I_-$ be the parts in any bipartition of the diagram.
For a simple root $\alpha_i$, let $\ep(i)$ be such that the simple generator corresponding to $\alpha_i$ is in $I_{\ep(i)}$.
For any $\alpha\in\Phi$, let $\alpha\ck$ be the corresponding {\em coroot} in the coroot system $\Phi\ck$.
The {\em fundamental weights} $\omega_i$ constitute the dual basis to the basis of simple coroots.  
That is, $\br{\alpha_i\ck,\omega_j}=\delta_{ij}$.
For any vector $x$ and any simple root $\alpha_i$, let $[x:\alpha_i]$ denote the coefficient of $\alpha_i$ when $x$ is expanded
in the basis $\Pi$.
Let $\Pge$ be the union of the positive roots and the negative simple roots of $\Phi$.
For $\ep\in\set{+,-}$, there is a map $\tau_\ep:\Pge\longrightarrow\Pge$ defined in~\cite{ga}.
Here we only note that every orbit of the cyclic group generated by $\tau_-\tau_+$ contains at least one negative simple root.
Thus for $\beta\in\Pge$, there is a well-defined {\em rotation number} $r(\beta)$ of $\beta$, defined to be the smallest integer $k$ such that 
$(\tau_-\tau_+)^k$ is a negative simple root.

The clusters are defined using $\tau_+$ and $\tau_-$.
A negative simple root $-\alpha_i$ is {\em compatible} with a root $\beta\in\Pge$ if and only if $[\beta:\alpha]=0$ and two positive roots 
$\beta$ and $\theta$ are compatible if and only if $\tau_\ep(\beta)$ and $\tau_\ep(\theta)$ are compatible for $\ep\in\set{+,-}$.
A {\em compatible subset} is a set of pairwise compatible elements of $\Pge$.
A {\em cluster} is a maximal compatible subset.
All clusters have cardinality $n$.
Call a set of $n-1$ pairwise compatible elements of $\Pge$ a {\em near-cluster}.
The {\em cluster fan} associated to $\Phi$ is the simplicial fan whose rays are the roots in $\Pge$ and whose maximal cones are generated by the 
clusters.
The maps $\tau_-$ and $\tau_+$ are combinatorial automorphisms of the cluster fan.
Chapoton, Fomin and Zelevinsky showed that the cluster fan is the normal fan of a polytope called the generalized associahedron for $\Phi$.

Suppose $C_1$ and $C_2$ are clusters with $|C_1\cap C_2|=n-1$, and define $\beta$ and $\theta$ so that $C_1-C_2=\set{\beta}$ and 
$C_2-C_1=\set{\theta}$.  
If $r(\beta)=r(\theta)=k$ then $(\tau_-\tau_+)^k\beta$ and $(\tau_-\tau_+)^k\theta$ are both negative simple roots, and thus they are 
compatible.
Therefore also $\beta$ and $\theta$ are compatible, so $C_1\cup C_2$ is a compatible subset of size $n+1$, which is impossible.
Therefore, without loss of generality we can take $r(\beta)<r(\theta)$ and define $C_1\le C_2$.
This defines an acyclic relation because whenever $C_1\covered C_2$ we have 
$\sum_{\alpha\in C_1}r(\alpha)<\sum_{\alpha\in C_2}r(\alpha)$.
Define the {\em cluster poset} to be the transitive closure of this relation.

Given a subset $\Pi'\subseteq\Pi$, the root subsystem $\Phi'$ of $\Phi$ generated by $\Pi'$ is called a {\em parabolic root subsystem}.
A pair of roots in $\Phi'_{\ge -1}$ is compatible in $\Phi'_{\ge -1}$ if and only if it is compatible in $\Pge$.

For vectors $x$ and $y$, define
\[\set{x,y}:=\sum_{\alpha_i\in\Pi}\ep(i)[x:\alpha_i\ck][y:\alpha_i].\]

The following lemma is proved in the course of the proof of~\cite[Proposition~3.1]{ga}.
\begin{lemma}
\label{twist}
If $\beta\ck$ is a positive coroot, $\theta$ is a positive root and $\ep\in\set{+,-}$, then
\[\set{\beta\ck,\tau_\ep\theta}=-\set{\tau_{-\ep}\beta\ck,\theta}.\]
\end{lemma}

The following lemma was previously observed independently by Zelevinsky~\cite{Zelev}.
To prove the lemma, it is useful to use some notation from~\cite{caII}.
Let $\tau_\ep^{(0)}$ be the identity map and let $\tau_\ep^{(k+1)}:=\tau_{(-1)^k\ep}\circ\tau_\ep^{(k)}$, so that for example 
$\tau_+^{(5)}=\tau_+\tau_-\tau_+\tau_-\tau_+$ and $\tau_+^{(6)}=\tau_-\tau_+\tau_-\tau_+\tau_-\tau_+$.

\begin{lemma}
\label{nice coroot}
For any near-cluster $N$, there is a positive coroot $\alpha\ck_N$ such that $\set{\alpha\ck_N,\alpha}=0$ for every $\alpha\in N$.
\end{lemma}
\begin{proof}
By induction on $n$.
If $n=1$, then $N=\emptyset$ and, setting $\alpha_N\ck=\alpha_1\ck$, the desired condition holds vacuously.
Suppose $n>1$, let $\beta$ and $\theta$ be the two distinct roots each of which completes $N$ to a cluster.
Let $k$ be the smallest natural number such that some root in $\tau_-^{(k)}(N\cup\set{\beta,\theta})$ is negative, and argue by induction on 
$k$ as well.
If $k=0$ and $-\alpha_i\in\set{\beta,\theta}$ for some $\alpha_i\in\Pi$ then the near-cluster $N$ is a subset of the parabolic root subsystem 
generated by $\Pi-\set{\alpha_i}$, so in particular $\set{\alpha_i\ck,\alpha}=0$ for every $\alpha\in N$.
If $k=0$ and $\beta,\theta\not\in-\Pi$ then there is some $-\alpha_i\in N$.
In this case $N':=N-\set{-\alpha_i}$ is a near-cluster in the parabolic root subsystem generated by $\Pi-\set{\alpha_i}$, and $\beta$ and 
$\theta$ are also contained in the subsystem.
Let the simple roots of the parabolic subsystem inherit the signing from the larger system.
By induction on $n$, there is a positive coroot $\alpha\ck_{N'}$ in the parabolic subsystem with 
$\set{\alpha\ck_{N'},\alpha}=0$ for every $\alpha\in N'$.
Since $\alpha\ck_{N'}$ is in the parabolic subsystem, we have $\set{\alpha\ck_{N'},-\alpha_i}=0$ as well.
Thus $\alpha_N\ck:=\alpha\ck_{N'}$ is the desired coroot.

If $k>0$ let $\ep=(-1)^k$ and consider $\tau_\ep N$.
By induction on $k$ there is a positive coroot $\alpha\ck_{(\tau_\ep N)}$ with 
$\set{\alpha\ck_{(\tau_\ep N)},\alpha}=0$ for every $\alpha\in \tau_\ep N$.
Since $\alpha\ck_{\tau_\ep N}$ is positive and each root in $N$ is positive, we apply Lemma~\ref{twist} to obtain 
$\set{\tau_{-\ep}\alpha\ck_{\tau_\ep N},\alpha}=0$ for every $\alpha\in N$.
If $\tau_{-\ep}\alpha\ck_{(\tau_\ep N)}$ is a negative simple coroot then $-\tau_{-\ep}\alpha\ck_{(\tau_\ep N)}$ is the desired positive root.
\end{proof}

The following simple linear-algebraic considerations are necessary for the proof of Theorem~\ref{cluster refine}.
Let $\set{b_i}$ be a basis for $\reals^n$, and let $\set{e_i}$ be the standard orthonormal basis.
Let $\br{\cdot,\cdot}$ be the standard inner product and let $\br{\cdot,\cdot}_b$ be the inner product for which $\set{b_i}$ is orthonormal,
that is $\br{x,y}_b=\br{Tx,Ty}$, where $T$ is the linear transformation mapping $b_i\mapsto e_i$ for all~$i$.
For $x\in\reals^n-\set{0}$, let $H(x)$ be the hyperplane normal to $x$ with respect to the standard inner product and let $H_b(x)$ be the 
hyperplane normal to $x$ with respect to $\br{\cdot,\cdot}_b$.
A very simple check shows that $H_b(x)=(T^*T)^{-1}H(x)$, where $T^*$ is the adjoint transformation to $T$
with respect to $\br{\cdot,\cdot}$.
\begin{comment}
\begin{eqnarray*}
H_b(x)&=&\set{y\in\reals^n:\br{x,y}_b=0}\\
&=&\set{y\in\reals^n:\br{Tx,Ty}=0}\\
&=&\set{y\in\reals^n:\br{x,T^*Ty}=0}\\
&=&\set{(T^*T)^{-1}y\in\reals^n:\br{x,y}=0}\\
&=&(T^*T)^{-1}\set{y\in\reals^n:\br{x,y}=0}\\
&=&(T^*T)^{-1}H(x)
\end{eqnarray*}
\end{comment}
Thus, given a collection $\Phi$ of vectors in $\reals^n$, the arrangement of hyperplanes normal to $\Phi$ with respect to $\br{\cdot,\cdot}_b$ is 
linearly isomorphic to the arrangement of hyperplanes normal to $\Phi$ with respect to the standard inner product.
Furthermore, given an isometry $\chi$ of $\br{\cdot,\cdot}_b$, the arrangement of hyperplanes normal to $\Phi$ with respect to 
$\br{\cdot,\cdot}_b$ is linearly isomorphic to the arrangement of hyperplanes normal to $\chi(\Phi)$ with respect to $\br{\cdot,\cdot}_b$.
The isometry is $\chi$ itself.

We now prove Theorem~\ref{cluster refine}, which states that the cluster fan for $W$ is refined by a fan which is linearly isomorphic to the 
normal fan of the $W$-permutohedron.
\begin{proof}[Proof of Theorem~\ref{cluster refine}]
The simple coroots $\set{\alpha_i\ck}$ are a basis of $\reals^n$, and thus the set $\set{c_i\alpha_i\ck}$ is a basis as well, where 
$\alpha_i=c_i^2\alpha_i\ck$.
Let $\br{\cdot,\cdot}_{c\alpha\ck}$ be the inner product for which $\set{c_i\alpha_i\ck}$ is an orthonormal basis.
For vectors $x$ and $y$ we can write
\begin{eqnarray*}
\br{x,y}_{c\alpha\ck}&=&\sum_{\alpha_i\in\Pi}[x:c_i\alpha_i\ck][y:c_i\alpha_i\ck]\\
&=&\sum_{\alpha_i\in\Pi}[x:\alpha_i\ck][y:c_i^2\alpha_i\ck]\\
&=&\sum_{\alpha_i\in\Pi}[x:\alpha_i\ck][y:\alpha_i]\\
\end{eqnarray*}
The normal fan of the $W$-permutohedron is defined by the arrangement $\A$ of hyperplanes normal, with respect to $\br{\cdot,\cdot}$, to the 
coroots $\Phi\ck$. 
This arrangement of hyperplanes is linearly isomorphic to the arrangement $\A'$ of hyperplanes normal, with respect to 
$\br{\cdot,\cdot}_{c\alpha\ck}$, to the coroots $\Phi\ck$. 

Define $\tw$ to be the linear ``twisting'' map which takes $\alpha_i\ck$ to $\ep(i)\alpha_i\ck$ for every~$i$.
Then $\tw$ is an isometry of $\br{\cdot,\cdot}_{c\alpha\ck}$, because it simply reverses the signs of some elements of an orthonormal 
basis.
Thus $\A'$ is in particular linearly isomorphic to the arrangement $\A''$ of hyperplanes normal, with respect to 
$\br{\cdot,\cdot}_{c\alpha\ck}$, 
to the twisted coroots $\tw(\Phi\ck)$. 
For $\beta\ck\in\Phi\ck$ and any vector $y$ we have
\begin{eqnarray*}
\br{\tw(\beta\ck),y}_{c\alpha\ck}&=&\sum_{\alpha_i\in\Pi}[\beta\ck:\ep(i)\alpha_i\ck][y:\alpha_i]\\
&=&\set{\beta\ck,y}
\end{eqnarray*}
Thus Lemma~\ref{nice coroot} says that every $(n-1)$-dimensional face of the cluster fan is contained in some hyperplane of $\A''$, so that
in particular, the cluster fan is refined by the fan defined by $\A''$, which is linearly isomorphic to the normal fan of the $W$-permutohedron.

The isometry is $T^*\circ T\circ\tw$, where $T$ is the map taking $c_i\alpha_i\ck$ to $e_i$, and $T^*$ is its adjoint.
One now easily checks that $T^*\circ T\circ\tw$ maps $\alpha_i$ to $\ep(i)\omega_i$.
\end{proof}

We conclude this section by relating the combinatorics of clusters for $S_n$ to the combinatorics of the Cambrian lattice for a bipartite 
orientation of the diagram of $S_n$.
Choose the linear functional $f$ and a polygon $Q$ so that each odd index is an up index and each even index is a down index.
This choice corresponds to the Cambrian lattice coming from the orientation $I_+\longrightarrow I_-$ for the signing $\ep(i)=(-1)^{i-1}$.
The map $\eta$ takes the identity permutation to a triangulation with diagonals $(1,2),(1,4),(3,4),(3,6),(5,6),(5,8)$ etc., where $(i,j)$ is 
shorthand for $(v_i,v_j)$.
In~\cite[Proposition 3.14]{ga}, the combinatorics of clusters for $S_n$ is made explicit with respect to the signing $\ep(i)=(-1)^{i-1}$.
As we quote this result, we relax the requirement of~\cite{ga} that $Q$ be regular and use the labeling of the vertices of $Q$ described 
earlier, rather than the labeling used in~\cite{ga}.
The roots in $\Pge$ are represented by the diagonals of $Q$ as follows.
The negative simple roots are the diagonals of the ``snake'' triangulation $\eta(12\cdots n)$.
Specifically, for~$i$ odd, $-\alpha_i$ corresponds to $(i,i+1)$ and for~$i$ even, $-\alpha_i$ is $(i-1,i+2)$.

We adopt the notation $\alpha_{i,j}:=\alpha_i+\cdots+\alpha_j$.
The labeling of the diagonals in the snake allows us to associate a positive root to each other diagonal $D$, by summing the negative simple
roots whose diagonals intersect $D$ in its interior, and then reversing sign.
Specifically, the positive root $\alpha_{i,j}$ corresponds to the diagonal
\[\left\lbrace\begin{array}{ll}
(i,j+1)&\mbox{if }i\mbox{ and }j\mbox{ are even,}\\
(i,j+3)&\mbox{if }i\mbox{ is even and }j\mbox{ is odd,}\\
(i-2,j+1)&\mbox{if }i\mbox{ is odd and }j\mbox{ is even, or}\\
(i-2,j+3)&\mbox{if }i\mbox{ and }j\mbox{ are odd.}\\
\end{array}\right.\]
However, if $i=1$, we replace $i-2$ by $0$ in the above correspondence and if $j=n-1$, we replace $j+3$ by $n+1$.
\begin{comment}
Specifically, if the diagonal $(i,j)$ is not in the snake, it corresponds to the root
\[\left\lbrace\begin{array}{ll}
\alpha_{i,j-3}&\mbox{if }i\mbox{ and }j\mbox{ are even,}\\
\alpha_{i,j-1}&\mbox{if }i\mbox{ is even and }j\mbox{ is odd,}\\
\alpha_{i+2,j-3}&\mbox{if }i\mbox{ is odd and }j\mbox{ is even, or}\\
\alpha_{i+2,j-1}&\mbox{if }i\mbox{ and }j\mbox{ are odd.}\\
\end{array}\right.\]
\end{comment}
The map $\tau_+$ is the combinatorial reflection of $Q$ which switches the endpoints of the edge $(0,2)$ and $\tau_-$ is the 
combinatorial reflection of $Q$ which switches the endpoints of the diagonal $(1,2)$.
The composition $\tau_-\tau_+$ is the combinatorial clockwise rotation of $Q$.
The clusters correspond to the sets of diagonals of triangulations of $Q$, so that two clusters intersect in a near-cluster if and only if they differ 
by a diagonal flip.
Each diagonal in the snake has a smaller slope than any other diagonal intersecting its interior, and combinatorial rotations preserve the 
comparison between the slopes of two diagonals unless one of them rotates through the vertical.
Thus the cluster poset is isomorphic to the Cambrian lattice $\C(I_+\longrightarrow I_-)$.

The following observations are useful in the next section.
For a positive root $\alpha_{i,j}$, the root $\tau_-\alpha_{i,j}$ corresponds to the diagonal $(i-1,j+2)$, regardless of the parity of 
$i$ and $j$, and no special statements need be made for the cases $i=1$ or $j=n-1$.
For a negative simple root, $-\alpha_i$, the root $\tau_-(-\alpha_i)=-\ep(i)\alpha_i$ corresponds to the diagonal $(i,i+1)$.

\section{Cambrian Fans}
\label{fans}
In this section we describe the Cambrian fans of types~A and~B in more detail.
For types~A and~B, we prove Conjecture~\ref{cluster conj} which states that the Cambrian fan for a bipartite orientation is linearly isomorphic 
to the cluster fan for the corresponding group.
We also show that the Cambrian fan associated to the Tamari lattice is the normal fan of the associahedron as realized by Stasheff~\cite{Sta95}.

For any fan $\F$, let the {\em rays} of $\F$ be the cones of $\F$ whose dimension is one larger than that of the minimal cone.
To make this definition match the usual definition of a ray, one can mod out by the minimal cone. 
Let $\A$ be a Coxeter arrangement for $W,$\ let $B$ be the region identified with the identity element of $W$ and let $-B$ be the antipodal 
region to $B$.
Given any region $R$ of $\A$, an {\em upper facet hyperplane} of $R$ is a facet hyperplane of $R$ which separates $R$ from a region above $R$
in the weak order.
The following proposition is a special case of~\cite[Proposition 5.10]{con_app}.
\begin{prop}
\label{rays}
Suppose that $\A$ is a Coxeter arrangement, that $\Theta$ does not contract any atoms and that $\F_\Theta$ is a simplicial fan.
Then the rays of $\F_\Theta$ are exactly the cones arising in one of the following two ways:
\begin{enumerate}
\item[(i) ]For a facet hyperplane $H$ of $B$, let $L$ be the subspace which is the intersection of the other facet hyperplanes of $B$.
Then the cone in $L$ consisting of points weakly separated from $-B$ by $H$ is a ray of $\F_\Theta$.
\item[(ii) ]Given a join-irreducible $\gamma$ of $\Po(\A,B)$ such that $\piup(\gamma)=\gamma$, let $L$ be the intersection of the upper 
facet hyperplanes of $\gamma$.
The cone consisting of points of~$L$ weakly separated from $B$ by the unique lower facet hyperplane of $\gamma$ is a ray of $\F(\Theta)$.
\end{enumerate}
\end{prop}

When $\A$ is the Coxeter arrangement for $S_n$ described in Section~\ref{weak}, these rays can be constructed explicitly.
To a ray of type (ii) we associate the subset $A$ corresponding to the join-irreducible $\gamma$.
To a ray of type (i) with $H$ normal to $e_{k+1}-e_k$ we associate the set $A=\set{k+1,\ldots,n}$.
Thus the rays of $\F$ correspond to proper non-empty subsets of $[n]$ and if $\Theta$ satisfies the hypotheses of Proposition~\ref{rays} then 
the rays of $\F_\Theta$ are the rays of $\F$ which do not correspond to a join-irreducible contracted by $\Theta$.
\begin{prop}
\label{explicit rays}
A non-empty proper subset $A\subset[n]$ corresponds to the ray with $x_i=x_j$ whenever $i,j\in A$ or $i,j\in A^c$ and $x_i<x_j$ whenever
$i\in A$ and $j\in A^c$.
\end{prop}
\begin{proof}
The upper facets of a join-irreducible $\gamma$ are the hyperplanes whose normals are $e_i-e_j$ for $i>j$ such that~$i$ immediately follows 
$j$ in the one-line notation for $\gamma$.
Thus the ray of type (ii) corresponding to $\gamma$ has $x_i=x_j$ whenever $i,j\in A$ or $i,j\in A^c$.
The unique lower facet of $\gamma$ is the hyperplane whose normal is $e_i-e_j$ for $i>j$ such that~$i$ is immediately followed by $j$ in 
$\gamma$.
We have $i\in A^c$ and $j\in A$.

If $A=[k+1,n]$ for some $k$ then the ray of type (i) corresponding to $A$ has $x_i=x_{i+1}$ for $i\in[n]-\set{k}$ and $x_k<x_{k+1}$.
\end{proof}

We now identify the join-irreducibles $\gamma$ with $\piup \gamma=\gamma$, where $\piup$ is the upward projection associated to 
a Cambrian congruence on $S_n$.
Recall from Proposition~\ref{projections} that a permutation $x$ has $\piup x=x$ if and only if $x$ avoids both $\up{2}13$ and $13\down{2}$.
For $1\le k\le n-1$, let $A_{k,k}$ be the set $[1,k]$ and for $1\le k<l\le n-1$, let $A_{k,l}$ be 
\[A_{k,l}:=\left\lbrace\begin{array}{ll}
[1,k]\cup\up{[k+1,l]}\cup[l+1,n]&\mbox{if }\down{k+1},\down{l}\\
\up{[k+1,l]}\cup[l+1,n]&\mbox{if }\up{k+1},\down{l}\mbox{\huge\protect\phantom{X}}\\
\ [1,k]\cup\up{[k+1,l]}&\mbox{if }\down{k+1},\up{l}\mbox{\huge\protect\phantom{X}}\\
\up{[k+1,l]}&\mbox{if }\up{k+1},\up{l}\mbox{\huge\protect\phantom{X}}
\end{array}\right.\]
\begin{prop}
\label{akl}
A join-irreducible avoids both $\up{2}13$ and $13\down{2}$ if and only if its corresponding subset is $A_{k,l}$ for some $1\le k\le l\le n-1$.
\end{prop}
\begin{proof}
Let $\gamma$ be a join-irreducible and let $A$ be its corresponding subset, with $m$ and $M$ as usual.
If $\gamma$ contains $a<\up{b}<c$ occurring in $\gamma$ in the order $bac$, then we may as well take $a$ to be $m$ and $c$ to be the maximum of $A$, so 
that $\up{b}$ is in $A^c\cap(m,\max A)$.
Conversely, if $A^c\cap(m,\max A)$ contains some $\up{b}$, then $\up{b}m(\max A)$ is a $\up{2}13$-pattern in $\gamma$.
Similarly, $\gamma$ contains a $13\down{2}$-pattern if and only if there is some $\down{b}\in A\cap(\min A^c,M)$.

Suppose $\gamma$ avoids both $\up{2}13$ and $13\down{2}$.
The following assertions are easily checked:  if $m>1$ then $\up{m}$; if $M<n$ then $\down{M}$; if $m=1$ and $\min A^c<\max A$ then 
$\down{\min A^c}$; and if $M=n$ and $\min A^c<\max A$ then $\up{\max A}$.
By the previous paragraph, we can construct $A$, breaking up into cases based on the identities of $m$, $M$, $\min A^c$ and $\max A$.
If $m>1$ and $M<n$ then $A=\up{[m,M]}\cup(M,n]$.
If $m=1$ and $M<n$ then $A=[1,\min A^c)\cup\up{[\min A^c,M]}\cup(M,n]$.
If $m>1$ and $M=n$ then $A=\up{[m,\max A]}$.
If $m=1$ and $M=n$ and $\min A^c<\max A$, then $A=[1,\min A^c)\cup\up{[\min A^c,\max A]}$.
Finally, if $m=1$ and $M=n$ and $\min A^c>\max A$, then $A=[1,\max A]$.
This last case is $A_{k,k}$ for $k=\max A$ and the first four cases correspond to the four cases defining $A_{k,l}$ for $k<l$.

Conversely, for $1\le k\le l\le n-1$ and $A=A_{k,l}$ one easily checks that \mbox{$A^c\cap\up{(m,\max A)}=\emptyset$} and 
$A\cap\down{(\min A^c,M)}=\emptyset$.
\end{proof}

Since a Cambrian fan $\F(\G)$ of type~A is combinatorially isomorphic to the normal fan of the associahedron, given any ray $F$ of $\F(\G)$,
we can find a diagonal of $Q$ such that every maximal cone containing $F$ corresponds to a triangulation using that diagonal.
The following proposition identifies the diagonal, and we refer to the diagonal $(v_i,v_j)$ simply as $(i,j)$. 
Given $i\in[n]$, let $\up{\nu}(i)$ be the smallest up index weakly greater than~$i$, and $\down{\nu}(i)$ the smallest down index 
weakly greater than~$i$, keeping in mind the fact that $0$ and $n+1$ are both up and down indices.
Similarly $\up{\mu}(i)$ and $\down{\mu}(i)$ are respectively the largest up and down indices weakly less than~$i$.
\begin{prop}
\label{diagonal}
The ray of type {\rm(i)} opposite the facet hyperplane $e_{k+1}-e_k$ of $B$ corresponds to the diagonal $(\up{\mu}(k),\down{\nu}(k+1))$.
The diagonal for a ray of type {\rm(ii)} corresponding to $A_{k,l}$ is $(\down{\mu}(k),\up{\nu}(k+1))$ if $k=l$ and if $k<l$ it is
\[\left\lbrace\begin{array}{ll}
(\down{\mu}(k),\down{\nu}(l+1)) &\mbox{if }\down{k+1},\down{l}\\
(\up{\mu}(k),\down{\nu}(l+1)) &\mbox{if }\up{k+1},\down{l}\mbox{\huge\protect\phantom{X}}\\
(\down{\mu}(k),\up{\nu}(l+1)) &\mbox{if }\down{k+1},\up{l}\mbox{\huge\protect\phantom{X}}\\
(\up{\mu}(k),\up{\nu}(l+1)) &\mbox{if }\up{k+1},\up{l}\mbox{\huge\protect\phantom{X}}
\end{array}\right.\]
\end{prop}
\begin{proof}
Recall from Section~\ref{A assoc} the description of the map $\eta$ from permutations to triangulations in terms of paths $\lambda$.
In the poset of regions corresponding to the weak order on $S_n$, a region contains the ray of type (i) opposite the facet hyperplane 
$e_{k+1}-e_k$ of $B$ if and only if every entry $\le k$ in the corresponding permutation occurs before every entry $\ge k+1$.
For every such permutation $x$, the path $\lambda_k(x)$ visits each up vertex with index $\le k$ then each down vertex with index $\ge k+1$.
Thus the diagonal $(\up{\mu}(k),\down{\nu}(k+1))$ occurs in every maximal cone of the Cambrian fan which contains this ray.

The rays of type (ii) correspond to the sets $A_{k,l}$ for $1\le k\le l\le n-1$.
Set $i:=n-|A_{k,l}|$.
In $S_n$, a region contains the ray corresponding to $A_{k,l}$ if and only if in the corresponding permutation $x$ every entry not in $A_{k,l}$
occurs before every entry in $A_{k,l}$.
If $l=k$ then for any such $x$, the path $\lambda_i(x)$ visits each down vertex in $[1,k]$ followed by each up vertex in $[k+1,n]$, and in
particular contains the diagonal $(\down{\mu}(k),\up{\nu}(k+1))$.
If $x$ is a permutation corresponding to a region containing the ray for $A_{k,l}$ with $k<l$, we check the four cases.
If $\down{k+1},\down{l}$ then the path $\lambda_i(x)$ visits each down vertex in $[1,k]$ followed by each down vertex in $[l+1,n]$, and thus 
contains the diagonal $(\down{\mu}(k),\down{\nu}(l+1))$.
If $\up{k+1},\down{l}$ then the path $\lambda_i(x)$ visits each up vertex in $[1,k]$ followed by each down vertex in $[l+1,n]$, and 
contains the diagonal $(\up{\mu}(k),\down{\nu}(l+1))$.
If $\down{k+1},\up{l}$ then the path $\lambda_i(x)$ visits each down vertex in $[1,k]$ followed by each up vertex in $[l+1,n]$, and 
contains the diagonal $(\down{\mu}(k),\up{\nu}(l+1))$.
If $\up{k+1},\up{l}$ then the path $\lambda_i(x)$ visits each up vertex in $[1,k]$ followed by each up vertex in $[l+1,n]$, and 
contains the diagonal $(\up{\mu}(k),\up{\nu}(l+1))$.
\end{proof}

We now explicitly determine the rays of the Cambrian fan of type~A which arises from one of the Tamari orientations.
Specifically, we take the congruence that arises from setting every index to be an up index.
Thus the join-irreducibles $\gamma$ with $\piup \gamma=\gamma$ are exactly the join-irreducibles arising from subsets either of the form 
$[1,k]$ for $1\le k\le n$, or $[k+1,l]$ for $1\le k<l\le n-1$.
The rays of type (i) correspond to subsets of the form $[k+1,n]$, so the rays of the fan correspond to proper subsets of the form $[i,j]$ for 
$1\le i\le j\le n$.
By Proposition~\ref{diagonal}, such a subset corresponds to the diagonal $(i-1,j+1)$ in $Q$.
Using Proposition~\ref{explicit rays} to explicitly write down the rays, we obtain an explicit description of the Cambrian fan for this Tamari 
orientation.
In Appendix B of~\cite{Sta95}, Stasheff and Shnider give a realization of the associahedron with these normal rays such that the correspondence
between rays and diagonals is the same.
Actually, their construction is phrased in terms of bracketings of a non-commutative product, but the well-known easy bijection between these
and triangulations of $Q$ gives the correct correspondence.
By the fact about Cambrian fans quoted in the introduction, the Tamari lattice is induced by a linear functional on this polytope.

We conclude by considering the Cambrian fan of type~A corresponding to a bipartite orientation and proving Conjecture~\ref{cluster conj} in 
types~A and~B.
Let $\ep(i)=(-1)^{i-1}$ and let $\G$ be the orientation $I_+\longrightarrow I_-$ of the diagram for $S_n$ defined at the end of 
Section~\ref{cluster}.
This orientation $\G$ corresponds to setting every odd index to be an up index and every even index to be a down index. 
We show that the Cambrian fan $\F(\G)$ is linearly isomorphic to the cluster fan, via the linear map discussed in the proof of 
Theorem~\ref{cluster refine}.
That is, we map the Cambrian fan to a refinement of the fan associated to the arrangement of hyperplanes normal to the twisted
root system $\tw(\Phi)$ for $S_n$, using the inner product $\br{\cdot,\cdot}_\alpha$ which makes the simple roots of 
$\Phi$ an orthonormal basis (recall that $\alpha\ck=\alpha$ for every root $\alpha$ in this root system).
We then show that this linearly transformed fan is the cluster fan.

We first find the rays of the transformed fan in the basis $\set{\alpha_i}$.
For each $k\in[n-1]$, there is a ray of type (i) consisting of points $x$ with $\br{x,\ep(k)\alpha_k}_\alpha\le 0$ which are orthogonal, with 
respect to $\br{\cdot,\cdot}_\alpha$, to $\ep(j)\alpha_j$ for all $j\neq k$.
In other words, the ray is generated by $-\ep(k)\alpha_k$.
The twisted upper normals of the join-irreducible corresponding to $A_{k,k}$ are $\set{\ep(j)\alpha_j:j\neq k}$, and the twisted lower normal is
$\tw(\alpha_{1,n-1})$, so the ray associated to $A_{k,k}$ is generated by $\ep(k)\alpha_k$.
For $k<l$ with $k$ odd and $l$ even, the twisted upper normals of the join-irreducible corresponding to $A_{k,l}$ are 
$\set{\ep(i)\alpha_i:i\in[1,k-1]\cup[l+1,n-1]}\cup\set{\tw(\alpha_{i,i+1}):i\in[k,l-1]}$ and the twisted lower normal is 
$\tw(\alpha_{1,l-1})$.
Thus the associated ray is $\alpha_{k,l}$.
A similar check of the other three cases finds the same the set of upper normals and the same associated ray, independent of the parity of $k$ or 
$l$.

We have shown that our linear transformation maps the rays of the Cambrian fan to the rays of the cluster fan, which are generated by twisted 
roots.
We use this correspondence to define a map $\psi$ from the diagonals of $Q$ to $\Pge$.
Specifically, we use Proposition~\ref{diagonal} to interpret the diagonal as a ray in the Cambrian fan and map that ray to the corresponding 
root in $\Pge$ by the linear transformation described above.
One checks that $\psi$ is the map $(a,a+1)\mapsto-\alpha_a$ and $(a,b)\mapsto\alpha_{a+1,b-2}$ if $b-a\ge 3$.
Comparing with the last paragraph of Section~\ref{cluster}, we see that $\psi$ is the inverse of the map which takes a root $\alpha$ to the 
diagonal associated to $\tau_-\alpha$.
Since $\tau_-$ is a combinatorial isomorphism of the cluster fan, we have shown that our linear transformation is also a combinatorial 
isomorphism.
Since the statement about the cluster poset was an easy observation at the end of Section~\ref{cluster}, this completes the proof of 
Conjecture~\ref{cluster conj} in the \mbox{type-A} case.

This result is easily extended to type~B.
The cluster fan for $B_n$ is the intersection of the cluster fan for $S_{2n}$ with the subspace of fixed points of the linear map $\chi$ which takes
$e_i\mapsto -e_{2n-i+1}$ for every $i\in[2n]$.
The map $\chi$ also takes $\alpha_i\mapsto\alpha_{2n-i}$ for every $i\in[2n-1]$, so its restriction to simple roots preserves the signing 
$\ep(i)=(-1)^{i-1}$.
In particular, $\chi$ commutes with the linear map which relates the cluster fan for $S_{2n}$ to the Cambrian fan for the bipartite orientation
of the diagram for $S_{2n}$.
Thus the Cambrian fan for the bipartite orientation of the diagram for $B_n$ is also obtained by intersecting with the subspace of fixed points
of $\chi$.

\section{Acknowledgments}
The author wishes to thank Sergey Fomin, Vic Reiner, John Stembridge and Hugh Thomas for helpful conversations.

\newcommand{\journalname}[1]{\textrm{#1}}
\newcommand{\booktitle}[1]{\textrm{#1}}

\end{document}